\documentclass{gtmon_a}
\pdfoutput=1
\usepackage[all]{xy}
\usepackage{pinlabel}

\proceedingstitle{Exotic homology manifolds (Oberwolfach 2003)}
\conferencestart{29 June 2003}
\conferenceend{5 July 2003}
\conferencename{Workshop on Exotic Homology Manifolds}
\conferencelocation{Oberwolfach Mathematics Institute, Oberwolfach,
Germany}

\editor{Frank Quinn}
\givenname{Frank}
\surname{Quinn}

\editor{Andrew Ranicki}
\givenname{Andrew}
\surname{Ranicki}

\title{The quadratic form $E_8$ and exotic homology manifolds}

\author{Washington Mio}
\givenname{Washington}
\surname{Mio}
\address{Department of Mathematics\\Florida State
University\\\newline
Tallahassee\\Florida 32306-4510\\USA}
\email{mio@math.fsu.edu}
\urladdr{}

\author{Andrew Ranicki}
\givenname{Andrew}
\surname{Ranicki}
\address{School of Mathematics\\University of Edinburgh\\\newline
Edinburgh\\EH9 3JZ\\Scotland\\United Kingdom}
\email{a.ranicki@ed.ac.uk}
\urladdr{}

\dedicatory{Dedicated to John Bryant on his 60th birthday}

\volumenumber{9}
\issuenumber{}
\publicationyear{2006}
\papernumber{4}
\startpage{33}
\endpage{66}

\doi{}
\MR{}
\Zbl{}

\arxivreference{math.GT/0403261}


\received{16 March 2004}
\revised{}
\accepted{16 March 2004}
\published{22 April 2006}
\publishedonline{22 April 2006}
\proposed{}
\seconded{}
\corresponding{}
\version{}

\makeatletter
\def\cnewtheorem#1[#2]#3{\newtheorem{#1}{#3}[section]
\expandafter\let\csname c@#1\endcsname\c@thm}

\newtheorem{thm}{Theorem}[section]
\cnewtheorem{lem}[thm]{Lemma}
\cnewtheorem{proposition}[thm]{Proposition}
\newtheorem*{thm81}{\fullref{T-E-eight}}
\theoremstyle{definition}
\cnewtheorem{defn}[thm]{Definition}
\cnewtheorem{rem}[thm]{Remark}
\cnewtheorem{example}[thm]{Example}

\makeatother  

\newcommand{\LL}{\ensuremath{\mathbb{L}}}

\newcommand{\forget}{\ensuremath{\mathscr{F}}}
\newcommand{\di}{\displaystyle}

\makeop{id}

\begin{document}

\begin{htmlabstract}
An explicit (-1)<sup>n</sup>&ndash;quadratic form over
<b>Z</b>[<b>Z</b><sup>2n</sup>] representing the surgery problem
E<sub>8</sub>&times;T<sup>2n</sup> is obtained, for use in the
Bryant&ndash;Ferry&ndash;Mio&ndash;Weinberger construction of
2n&ndash;dimensional exotic homology manifolds.
\end{htmlabstract}

\begin{abstract} 
An explicit $(-1)^n$--quadratic form over $\mathbb{Z}[\mathbb{Z}^{2n}]$
representing the surgery problem $E_8 \times T^{2n}$ is obtained,
for use in the Bryant--Ferry--Mio--Weinberger construction of
$2n$--dimensional exotic homology manifolds.
\end{abstract}

\maketitle

%

\section{Introduction} \label{S-intro}

Exotic ENR homology $n$--manifolds, $n \geqslant 6$,
were discovered in the early 1990s by Bryant, Ferry, Mio and Weinberger
\cite{BFMW93,BFMW96}.  In the 1970s, the existence of such spaces
had become a widely debated problem among geometric topologists in
connection with the works of Cannon \cite{jC79}, Edwards \cite{rE80} and
Daverman \cite{rD86} on the characterization
of topological manifolds.
The {\em Resolution Conjecture}, formulated by Cannon in \cite{jC78},
implied the non-existence of exotic homology manifolds --
compelling evidence supporting the conjecture was offered by  the
solution of the {\em Double Suspension Problem}.  Quinn introduced
methods of controlled $K$--theory and controlled surgery into the area.
He associated with an ENR
homology $n$--manifold $X$, $n \geqslant 5$, a local index
$\imath(X) \in 8 \mathbb{Z} + 1$ with the property that $\imath(X) =1$
if and only if $X$ is resolvable.  A {\em resolution\,} of $X$ is a
proper surjection $f \colon M \to X$ from a topological manifold $M$
such that, for each $x \in X$, $f^{-1} (x)$ is contractible in any of its
neighborhoods in $M$.  This led to the celebrated Edwards--Quinn
characterization of topological $n$--manifolds, $n \geqslant 5$, as index--1
ENR homology manifolds satisfying the disjoint disks property (DDP)
(see the articles by Quinn \cite{fQ83,fQ87} and Daverman \cite{rD86}.
More details and historical remarks on these
developments can be found in the survey articles by Cannon \cite{jC78},
Edwards \cite{rE80}, Weinberger \cite{sW95} and Mio \cite{wM00}, and in
the book by Daverman \cite{rD86}.

In \cite{BFMW93,BFMW96}, ENR homology manifolds with non-trivial local indexes
are constructed as inverse limits of ever finer Poincar\'e duality
spaces, which are obtained from topological manifolds using controlled
cut-paste constructions.  In the simply-connected case, for example,
topological manifolds are cut along the boundaries of regular
neighborhoods of very fine 2--skeleta and pasted back together using
$\epsilon$--homotopy equivalences that ``carry non-trivial local
indexes'' in the form of obstructions to deform them to homeomorphisms
in a controlled manner.  The construction of these
$\epsilon$--equivalences requires controlled surgery theory, the
calculation of controlled surgery groups with trivial local fundamental
group, and ``Wall realization'' of controlled surgery obstructions.
The stability of controlled surgery groups is a key fact, whose proof
was completed more recently by Pedersen, Quinn and Ranicki
\cite{PQR-zero-one}; an elegant proof along similar lines was given by
Pedersen and Yamasaki \cite{PY} at the 2003 Workshop on Exotic Homology
Manifolds in Oberwolfach, employing methods of Yamasaki \cite{mY87}.  An
alternative proof based on the {\em $\alpha$--Approximation Theorem\/}
is due to Ferry \cite{sF-zero-three}.

The construction of exotic homology manifolds presented by Bryant, Ferry,
Mio and Weinberger \cite{BFMW96} is somewhat indirect.  Along the years,
many colleagues
(notably Bob Edwards) voiced the desire to see -- at least in one
specific example -- an explicit realization of the controlled quadratic
form employed in the Wall realization of the local index.  This became
even clearer at the workshop in Oberwolfach.  A detailed inspection of
the construction of \cite{BFMW96} reveals that it suffices to give this
explicit description at the first (controlled) stage of the
construction of the inverse limit, since fairly general arguments show
that subsequent stages can be designed to inherit the local index.

The
main goal of this paper is to provide explicit realizations of
controlled quadratic forms that lead to the construction of compact
exotic homology manifolds with fundamental group $\Z^{2n}$, $n
\geqslant 3$, which are not homotopy equivalent to any closed
topological manifold.  This construction was suggested in 
\cite[Section~7]{BFMW96}, but details were not provided.  Starting with
the rank--8
quadratic form $E_8$ of signature 8, which generates the Wall group
$L_0 (\Z) \cong \Z$, we explicitly realize its image in $L_{2n}
(\Z[\Z^{2n}])$ under the canonical embedding $L_0 (\Z) \to L_{2n}
(\Z[\Z^{2n}])$.

Let
$$\psi_0=\begin{pmatrix}
1&0&0&1&0&0&0&0 \\
0&1&1&0&0&0&0&0 \\
0&0&1&1&0&0&0&0 \\
0&0&0&1&1&0&0&0 \\
0&0&0&0&1&1&0&0 \\
0&0&0&0&0&1&1&0 \\
0&0&0&0&0&0&1&1 \\
0&0&0&0&0&0&0&1
\end{pmatrix}$$
be the $8 \times 8$ matrix over $\Z$ with symmetrization
the unimodular $8 \times 8$ matrix of the $E_8$--form:
$$\psi_0+\psi_0^*=E_8=\begin{pmatrix}
2&0&0&1&0&0&0&0 \\
0&2&1&0&0&0&0&0 \\
0&1&2&1&0&0&0&0 \\
1&0&1&2&1&0&0&0 \\
0&0&0&1&2&1&0&0 \\
0&0&0&0&1&2&1&0 \\
0&0&0&0&0&1&2&1 \\
0&0&0&0&0&0&1&2
\end{pmatrix}$$
Write
$$\begin{array}{ll}
\Z[\Z^{2n}]&=\Z[z_1,z_1^{-1},\dots,z_{2n},z_{2n}^{-1}]\\[1ex]
&=\Z[z_1,z_1^{-1}]\otimes \Z[z_2,z^{-1}_2] \otimes
\dots \otimes \Z[z_{2n},z_{2n}^{-1}].
\end{array}$$
For $i=1,2,\dots,n$ define the $2 \times 2$ matrix over
$\Z[z_{2i-1},z^{-1}_{2i-1},z_{2i},z^{-1}_{2i}]$
$$\alpha_i=\begin{pmatrix} 1-z_{2i-1} & z_{2i-1}z_{2i}-z_{2i-1}-z_{2i} \\
1 & 1-z_{2i} \end{pmatrix},$$
so that $\alpha_1 \otimes \alpha_2 \otimes \dots \otimes \alpha_n$ is
a $2^n \times 2^n$ matrix over $\Z[\Z^{2n}]$. (See \fullref{S-almost}
for the geometric provenance of the matrices $\alpha_i$).

\begin{thm81}
The surgery obstruction
$E_8\times T^{2n} \in L_{2n}(\Lambda)$ $(\Lambda=\Z[\Z^{2n}])$
is represented by the nonsingular
$(-1)^n$--quadratic form $(K,\lambda,\mu)$ over $\Lambda$,
with
$$K=\Z^8 \otimes \Lambda^{2^n}=\Lambda^{2^{n+3}}$$
the f.g.\ free $\Lambda$--module of rank $8.2^n=2^{n+3}$ and
$$\begin{array}{l}
\lambda=\psi+(-1)^n\psi^*\co K \to K^*=
{\rm Hom}_{\Lambda}(K,\Lambda),\\[1ex]
\mu(x)=\psi(x)(x) \in Q_{(-1)^n}(\Lambda)=
\Lambda/\{g+(-1)^{n+1}g^{-1}\,\vert\,g \in \Z^{2n}\}~(x\in K)
\end{array}
$$
with
$$\psi=\psi_0\otimes \alpha_1 \otimes \alpha_2 \otimes
\dots \otimes \alpha_n\co K\to  K^*.$$
\end{thm81}

Sections \ref{S-wall}--\ref{S-explicit} contain background material on
surgery theory
and the arguments that lead to a proof of \fullref{T-E-eight}.
Invariance of $E_8 \times T^{2n}$ under transfers to finite covers
is proven in \fullref{S-transfer}.
In \fullref{S-control}, using a large finite cover $T^{2n} \to T^{2n}$,
we describe how to pass from the
non-simply-connected surgery obstruction $E_8 \times T^{2n}$  to
a controlled quadratic
$\Z$--form over $T^{2n}$.  Finally, in \fullref{S-exotic} we explain how
the controlled version of $E_8 \times T^{2n}$ is used in the construction
of exotic homology $2n$--manifolds $X$ with Quinn index $\imath(X)=9$.

\section{The Wall groups} \label{S-wall}

We begin with some recollections of surgery obstruction theory --
we only need the details in the even-dimensional oriented case.

Let $\Lambda$ be a ring with an involution, that is a function
$$\raise5pt\hbox{$\overline{\hphantom{a}}$}\co
\Lambda \to \Lambda;\ a \mapsto \overline{a}$$
satisfying
$$\overline{a+b}=\overline{a}+\overline{b},\quad
\overline{ab}=\overline{b}.\overline{a},\quad
\overline{\overline{a}}=a,\quad\overline{1}=1 \in \Lambda.$$

\begin{example}
In the applications to topology $\Lambda=\Z[\pi]$ is a group ring,
with the involution
$$\raise5pt\hbox{$\overline{\hphantom{a}}$}\co
\Z[\pi] \to  \Z[\pi];\quad\sum\limits_{g \in \pi}a_gg \mapsto
\sum\limits_{g \in \pi}a_gg^{-1}\quad (a_g \in \Z).$$
\end{example}

The involution is used to define a left $\Lambda$--module structure on
the dual of a left $\Lambda$--module $K$
$$K^*:={\rm Hom}_{\Lambda}(K,\Lambda),$$
with
$$\Lambda \times K^* \to K^*;\ (a,f) \mapsto (x \mapsto f(x).\overline{a}).$$

The {\it $2n$--dimensional surgery obstruction group} $L_{2n}(\Lambda)$
is defined by Wall \cite[Section~5]{cW-seventy} to be the Witt
group of {\it nonsingular $(-1)^n$--quadratic forms} $(K,\lambda,\mu)$
over $\Lambda$, with $K$ a finitely generated free (left)
$\Lambda$--module together with
\begin{itemize}
\item[(i)] a pairing
$\lambda \colon K \times K \to \Lambda$
such that
\begin{align*}
\lambda(x,ay)&=a\lambda(x,y),\\
\lambda(x,y+z)&=\lambda(x,y)+\lambda(x,z),\\
\lambda(y,x)&=(-1)^n\overline{\lambda(x,y)}
\end{align*}
and the adjoint $\Lambda$--module morphism
$$\lambda\co K \to K^*;\quad x \mapsto (y \mapsto \lambda(x,y))$$
is an isomorphism,
\item[(ii)] a $(-1)^n$--quadratic function
$\mu \colon K \to Q_{(-1)^n}(\Lambda) = \Lambda/
\{a +(-1)^{n+1}\overline{a}\,\vert\, a \in \Lambda\}$
with
\begin{align*}
\lambda(x,x)&=\mu(x)+(-1)^n\overline{\mu(x)},\\
\mu(x+y)&=\mu(x)+\mu(y)+\lambda(x,y),\\
\mu(ax)&=a\mu(x)\overline{a}.
\end{align*}
\end{itemize}
For a f.g.\ free $\Lambda$--module $K=\Lambda^r$
with basis $\{e_1,e_2,\dots,e_r\}$ the pair $(\lambda,\mu)$
can be regarded as an equivalence class of $r \times r$ matrices over $\Lambda$
\[
\psi = (\psi_{ij})_{1 \leqslant i,j \leqslant r}\quad (\psi_{ij} \in \Lambda)
\]
such that $\psi+(-1)^n\psi^*$ is invertible, with $\psi^*=(\overline{\psi}_{ji})$,
and
$$\psi \sim \psi'~{\rm if}~\psi'-\psi = \chi+(-1)^{n+1}\chi^*~{\rm for~some}~
r \times r~{\rm matrix}~\chi=(\chi_{ij}).$$
The relationship between $(\lambda,\mu)$ and $\psi$ is given by
$$\begin{array}{l}
\lambda(e_i,e_j) = \psi_{ij}+(-1)^n\overline{\psi}_{ji} \in \Lambda,\\[1ex]
\mu(e_i) = \psi_{ii} \in Q_{(-1)^n}(\Lambda),
\end{array}$$
and we shall write
$$(\lambda,\mu) = (\psi+(-)^n\psi^*,\psi).$$
\indent
The detailed definitions of the odd-dimensional $L$--groups $L_{2n+1}(\Lambda)$
are rather more complicated, and are not required here.
The quadratic $L$--groups are 4--periodic
$$L_m(\Lambda) = L_{m+4}(\Lambda).$$
\indent The simply-connected quadratic $L$--groups are given by
$$L_m(\Z) = P_m = \begin{cases}
\Z&{\rm if}~m \equiv 0(\bmod\,4)\\
0&{\rm if}~m \equiv 1(\bmod\,4)\\
\Z_2&{\rm if}~m \equiv 2(\bmod\,4)\\
0&{\rm if}~m \equiv 3(\bmod\,4)
\end{cases}$$
(Kervaire--Milnor). In particular, for $m\equiv 0(\bmod\,4)$ there
is defined an isomorphism
$$L_0(\Z) \xymatrix{\ar[r]^-{\di{\cong}}&} \Z;\quad (K,\lambda,\mu)
\mapsto \tfrac{1}{8}\,{\rm signature}(K,\lambda).$$
\indent
The {\it kernel form} of an $n$--connected normal map
$(f,b)\co M^{2n} \to X$ from a $2n$--dimensional manifold $M$ to an
oriented $2n$--dimensional geometric Poincar\'e complex $X$
is the nonsingular $(-1)^n$--quadratic form
over $\Z[\pi_1(X)]$ defined by Wall \cite[Section~5]{cW-seventy}
$$(K_n(M),\lambda,\mu)$$
with
$$K_n(M) = {\rm ker}(\wtilde{f}_*\co H_n(\wwtilde{M}) \to H_n(\wwtilde{X}))$$
the kernel (stably) f.g.\ free $\Z[\pi_1(X)]$--module, $\wwtilde{X}$
the universal cover of $X$, $\wwtilde{M}=f^*\wwtilde{X}$ the pullback cover
of $M$ and $(\lambda,\mu)$ given by geometric (intersection, self-intersection) numbers.
The {\it surgery obstruction} of Wall \cite{cW-seventy}
$$\sigma_*(f,b) = (K_n(M),\lambda,\mu) \in L_{2n}(\Z[\pi_1(X)])$$
is such that $\sigma_*(f,b)=0$ if (and for $n \geqslant 3$ only if)
$(f,b)$ is bordant to a homotopy equivalence.

The Realization Theorem \cite[Section~5]{cW-seventy} states that for a
finitely presented group $\pi$ and $n\geqslant 3$ every nonsingular
$(-1)^n$--quadratic form $(K,\lambda,\mu)$ over $\Z[\pi]$ is the
kernel form of an $n$--connected $2n$--dimensional normal map
$f\co M \to X$ with $\pi_1(X)=\pi$.

\section{The instant surgery obstruction}\label{S-instant}

Let $(f,b)\co M^m \to X$ be an $m$--dimensional normal map with
$f_*\co \pi_1(M) \to \pi_1(X)$ an isomorphism, and let
$\wtilde{f}\co \wwtilde{M}\to \wwtilde{X}$ be a
$\pi_1(X)$--equivariant lift of $f$ to the universal covers of
$M,X$. The $\Z[\pi_1(X)]$--module morphisms
$\wtilde{f}_*\co H_r(\wwtilde{M}) \to H_r(\wwtilde{X})$ are
split surjections, with the Umkehr maps
$$f^!\co H_r(\wwtilde{X}) \cong H^{m-r}(\wwtilde{X})
\xymatrix{\ar[r]^-{\di{\wtilde{f}^*}}&}
    H^{m-r}(\wwtilde{M}) \cong H_r(\wwtilde{M}),$$
such that
$$\wtilde{f}_*f^!=1\co H_r(\wwtilde{X}) \to
    H_r(\wwtilde{X}).$$
The kernel $\Z[\pi_1(X)]$--modules
$$K_r(M)={\rm ker}(\wtilde{f}_*\co
  H_r(\wwtilde{M}) \to H_r(\wwtilde{X}))$$
are such that
$$H_r(\wwtilde{M})=K_r(M) \oplus
  H_r(\wwtilde{X}),K_r(M)=H_{r+1}(\wtilde{f}).$$
By the Hurewicz theorem, $(f,b)$ is $k$--connected if and only if
$$K_r(M) = 0~{\rm for}~r<k,$$
in which case $K_k(M)=\pi_{k+1}(f)$. If $m=2n$ or $2n+1$ then by
Poincar\'e duality $(f,b)$ is $(n+1)$--connected if and only if it
is a homotopy equivalence. In the even-dimensional case $m=2n$ the
surgery obstruction of $(f,b)$ is defined to be
$$\sigma_*(f,b) = \sigma_*(f',b') = (K_n(M'),\lambda',\mu') \in L_{2n}(\Z[\pi_1(X)])$$
with $(f',b'):M' \to X$ any bordant
$n$--connected normal map obtained from $(f,b)$ by surgery below
the middle dimension.  The instant surgery obstruction of
Ranicki~\cite{ranicki-one} is an expression for such a form
$(K_n(M'),\lambda',\mu')$ in terms of the kernel $2n$--dimensional quadratic
Poincar\'e complex $(C,\psi)$ such that $H_*(C)=K_*(M)$.  In
\fullref{S-explicit} we below we shall use a variant of the instant
surgery obstruction to obtain an explicit $(-1)^n$--quadratic form
over $\Z[\Z^{2n}]$ representing $E_8\times T^{2n}
\in L_{2n}(\Z[\Z^{2n}])$.

Given a ring with involution $\Lambda$ and an $m$--dimensional
f.g.\ free $\Lambda$--module chain complex
$$C\co C_m \xymatrix{\ar[r]^-{\di{d}}&} C_{m-1} \to \dots \to
C_1 \xymatrix{\ar[r]^-{\di{d}}&} C_0$$
let $C^{m-*}$ be the dual $m$--dimensional f.g.\ free $\Lambda$--module chain complex,
with
$$\begin{array}{l}
d_{C^{m-*}} = (-1)^rd^*\co \\[1ex]
(C^{m-*})_r = C^{m-r} = (C_{m-r})^* = 
{\rm Hom}_{\Lambda}(C_{m-r},\Lambda)\to C^{m-r+1}.
\end{array}$$
Define a duality involution on ${\rm Hom}_{\Lambda}(C^{m-*},C)$ by
$$T\co {\rm Hom}_{\Lambda}(C^p,C_q) \to {\rm Hom}_{\Lambda}(C^q,C_p);\quad 
\phi \mapsto (-1)^{pq}\phi^*.$$
An {\it $m$--dimensional quadratic Poincar\'e complex} $(C,\psi)$ over $\Lambda$
is an $m$--dimensional f.g.\ free $\Lambda$--module chain complex $C$
together with $\Lambda$--module morphisms
$$\psi_s\co C^r \to C_{m-r-s}~~
(s \geqslant 0)$$
such that
$$d\psi_s+(-1)^r\psi_sd^*+(-1)^{m-s-1}(\psi_{s+1}+(-1)^{s+1}T\psi_{s+1}) = 0\co 
C^{m-r-s-1} \to C_r$$
for $s\geqslant 0$,
and such that $(1+T)\psi_0\co C^{m-*} \to C$ is a chain equivalence.
The cobordism group of $m$--dimensional quadratic Poincar\'e
complexes over $\Lambda$ was identified in Ranicki \cite{ranicki-one}
with the Wall surgery obstruction $L_m(\Lambda)$, and the surgery obstruction
of an $m$--dimensional normal map $(f,b)\co M \to X$
was identified with the cobordism class
$$\sigma_*(f,b) = ({\mathcal C}(f^!),\psi_b) \in L_m(\Z[\pi_1(X)])$$
of the kernel quadratic Poincar\'e complex $({\mathcal C}(f^!),\psi_b)$,
with ${\mathcal C}(f^!)$ the algebraic mapping cone of the Umkehr
$\Z[\pi_1(X)]$--module chain map
$$f^!\co C(\wwtilde{X}) \simeq C(\wwtilde{X})^{m-*}
\xymatrix{\ar[r]^-{\di{\wtilde{f}^*}}&} C(\wwtilde{M})^{m-*}~
\simeq~ C(\wwtilde{M}).$$
The homology $\Z[\pi_1(X)]$--modules of ${\mathcal C}(f^!)$ are
the kernel $\Z[\pi_1(X)]$--modules of $f$
$$H_*({\mathcal C}(f^!)) = K_*(M) = {\rm ker}(\wtilde{f}_*\co H_*(\wwtilde{M})
\to H_*(\wwtilde{X})).$$

\begin{defn} The {\it instant form} of a $2n$--dimensional
quadratic Poincar\'e complex $(C,\psi)$ over $\Lambda$
is the nonsingular $(-1)^n$--quadratic form over $\Lambda$
$$\begin{array}{l}
(K,\lambda,\mu) = 
\bigg({\rm coker}(\begin{pmatrix}
d^* & 0 \\ (-1)^{n+1}(1+T)\psi_0 & d \end{pmatrix}\co 
C^{n-1} \oplus C_{n+2}  \to C^n \oplus C_{n+1}),\\[2ex]
\hskip175pt
\bigg[\begin{matrix} \psi_0+(-1)^n\psi^*_0 & d \\
(-1)^nd^* & 0 \end{matrix}\bigg],
\bigg[\begin{matrix} \psi_0 & d \\
0 & 0 \end{matrix}\bigg]\bigg).
\end{array}$$
\end{defn}

If $C_r$ is f.g.\ free with ${\rm rank}_{\Lambda}C_r=c_r$ then $K$ is
(stably) f.g.\ free with
$${\rm rank}_{\Lambda}K = \sum\limits^n_{r=0}(-1)^r(c_{n-r}+c_{n+r+1})\in \Z.$$
If $(1+T)\psi_0\co C^{2n-*} \to C$ is an isomorphism then
$$c_{n+r+1}=c_{n-r-1},~{\rm rank}_{\Lambda}K = c_n,$$
with
$$(K,\lambda,\mu) = (C^n,\psi_0+(-1)^n\psi^*_0,\psi_0).$$

\begin{proposition}[Instant surgery obstruction, Ranicki
{{\cite[Proposition I.4.3]{ranicki-one}}}]\qua\newline
{\rm (i)}\qua The cobordism class of a $2n$--dimensional
quadratic Poincar\'e complex $(C,\psi)$ over $\Lambda$
is the Witt class
$$(C,\psi) = (K,\lambda,\mu) \in L_{2n}(\Lambda)$$
of the instant nonsingular $(-1)^n$--quadratic form $(K,\lambda,\mu)$
over $\Lambda$.

{\rm (ii)}\qua The surgery obstruction of a $2n$--dimensional normal map
$(f,b)\co M \to X$ is represented by the instant form $(K,\lambda,\mu)$ of any
quadratic Poincar\'e complex $(C,\psi)$ which is chain equivalent to
the kernel
$2n$--dimensional quadratic  Poincar\'e complex $(C(f^!),\psi_b)$
$$\sigma_*(f,b) = (K,\lambda,\mu) \in L_{2n}(\Z[\pi_1(X)]).$$
\end{proposition}

\begin{rem}  \label{below}
(i)\qua If $(f,b)$ is $n$--connected then $C$ is chain equivalent to the chain
complex concentrated in dimension $n$
$$C\co 0 \to \dots \to 0 \to K_n(M) \to 0 \to  \dots \to 0$$
and the instant form is just the kernel form $(K_n(M),\lambda,\mu)$
of Wall \cite{cW-seventy}.

(ii)\qua More generally, if
$(f,b)$ is $k$--connected for some $k \leqslant n$
then $C$ is chain equivalent to a chain complex concentrated
in dimensions $k,k+1,\dots,2n-k$
$$C\co 0 \to \dots \to 0 \to C_{2n-k} \to \dots \to
C_n \to \dots \to C_k \to 0 \to  \dots \to 0.$$
For $n \geqslant 3$ the effect of surgeries killing the
$c_{2n-k}$ generators of $H^{2n-k}(C)=K_k(M)$ represented by a basis
of $C^{2n-k}$ is a bordant $(k+1)$--connected normal map
$$(f',b')\co M'\to X$$
with ${\mathscr C}(f^{\prime !}\co C(\wwtilde{X})
\to C(\wwtilde{M}'))$ chain equivalent to a chain complex of the type
$$C'\co 0 \to \dots \to 0 \to  C'_{2n-k-1} \to \dots \to C'_n \to \dots
\to C'_{k+1} \to 0 \to \dots \to 0$$
with
$$C'_r = 
\begin{cases}
{\rm ker}((d~(1+T)\psi_0)\co C_{k+1} \oplus C^{2n-k} \to C_k)&{\rm if}~r=k+1\\[1ex]
C_r&{\rm if}~k+2 \leqslant r \leqslant 2n-k-1.
\end{cases}$$
Proceeding in this way, there is obtained a sequence of bordant
$j$--connected normal maps
$$(f_j,b_j)\co M_j \to X~~(j=k,k+1,\dots,n)$$
with
$$(f_k,b_k) = (f,b),~(f_{j+1},b_{j+1}) = (f_j,b_j)'.$$
The instant form of $(C,\psi)$ is precisely
the kernel $(-1)^n$--quadratic form
$$(K_n(M_n),\lambda_n,\mu_n)$$
of the $n$--connected normal map $(f_n,b_n)\co M_n \to X$, so that
the surgery obstruction of $(f,b)$ is given by
$$\begin{array}{ll}
\sigma_*(f,b)&= \sigma_*(f_k,b_k) = \dots = \sigma_*(f_n,b_n)\\[1ex]
&= (K_n(M_n),\lambda_n,\mu_n) \in L_{2n}(\Z[\pi_1(X)]).
\end{array}$$
\end{rem}

\section{The quadratic form $E_8$}

For $m \geqslant 2$ let $M_0^{4m}$ be the $(2m-1)$--connected
$4m$--dimensional $PL$ manifold obtained from the Milnor $E_8$--plumbing
of 8 copies of $\tau_{S^{2m}}$ by coning off the (exotic)
$(4m-1)$--sphere boundary, with intersection form $E_8$ of signature 8.
(For $m=1$ we can take $M_0$ to be the simply-connected 4--dimensional Freedman
topological manifold with intersection form $E_8$).  The surgery
obstruction of the corresponding $2m$--connected normal map
$(f_0,b_0)\co M_0^{4m} \to S^{4m}$ represents the generator
$$\sigma_*(f_0,b_0) = (K_{2m}(M_0),\lambda,\mu) = (\Z^8,E_8) = 1 \in L_{4m}(\Z) = L_0(\Z) = \Z$$
with
$$\begin{array}{l}
K_{2m}(M_0) = H_{2m}(M_0) = \Z^8\\[2ex]
\lambda = E_8 = \begin{pmatrix}
2&0&0&1&0&0&0&0 \\
0&2&1&0&0&0&0&0 \\
0&1&2&1&0&0&0&0 \\
1&0&1&2&1&0&0&0 \\
0&0&0&1&2&1&0&0 \\
0&0&0&0&1&2&1&0 \\
0&0&0&0&0&1&2&1 \\
0&0&0&0&0&0&1&2
\end{pmatrix}\\[11ex]
\mu(0,\dots,1,\dots,0) = 1.
\end{array}$$

\section{The surgery product formula}\label{S-product}

Surgery product formulae were originally obtained in the
simply-connected case, notably by Sullivan.
We now recall the non-simply-connected
surgery product formula of Ranicki \cite{ranicki-one} involving the
Mishchenko symmetric $L$--groups. In \fullref{S-almost} we shall
recall the variant of the surgery product formula involving
almost symmetric $L$--groups of Clauwens,
which will be used in \fullref{T-E-eight} below to write down an
explicit nonsingular $(-1)^n$--quadratic form over $\Z[\Z^{2n}]$ $(n
\geqslant 1)$ representing the image of the generator
\[
1 = E_8 \in L_{4m}(\Z) \cong \Z~(m \geqslant 0)
\]
under the canonical embedding
\begin{multline*}
{}_{-}\times T^{2n} \colon L_{4m}(\Z) \to L_{4m+2n}(\Z[\Z^{2n}]);\quad \\[1ex]
\sigma_*((f_0,b_0)\co M_0\to S^{4m})=E_8 \mapsto
\sigma_*((f_0,b_0)\times 1\co M_0\times T^{2n} \to S^{4m}\times T^{2n}).
\end{multline*}

An {\it $n$--dimensional symmetric Poincar\'e complex}
$(C,\phi)$ over a ring with involution $\Lambda$
is an $n$--dimensional f.g.\ free $\Lambda$--module chain complex
$$C\co C_{n} \xymatrix{\ar[r]^-{\di{d}}&
C_{n-1} \ar[r]& \dots \ar[r] & C_1 \ar[r]^-{\di{d}}&C_0}$$
together with $\Lambda$--module morphisms
$$\phi_s\co C^r = {\rm Hom}_{\Lambda}(C_r,\Lambda) \to C_{n-r+s}~~
(s \geqslant 0)$$
such that
$$\begin{array}{l}
d\phi_s+(-1)^r\phi_sd^*+(-1)^{n+s-1}(\phi_{s-1}+(-1)^sT\phi_{s-1})=0\co\\[1ex]
\hskip170pt C^{n-r+s-1} \to C_r~~ (s\geqslant 0,\phi_{-1}=0)
\end{array}$$
and $\phi_0\co C^{n-*} \to C$ is a chain equivalence.
The cobordism group of $n$--dimensional symmetric
Poincar\'e complexes over $\Lambda$ is denoted by $L^n(\Lambda)$
-- see Ranicki \cite{ranicki-one} for a detailed exposition of symmetric
$L$--theory.
Note that the symmetric $L$--groups $L^*(\Lambda)$ are not 4--periodic
in general
$$L^n(\Lambda) \neq L^{n+4}(\Lambda).$$
\indent The symmetric $L$--groups of $\Z$ are given by
$$L^n(\Z) = \begin{cases}
\Z&{\rm if}~n \equiv 0(\bmod\,4)\\
\Z_2&{\rm if}~n \equiv 1(\bmod\,4)\\
0&{\rm if}~n \equiv 2(\bmod\,4)\\
0&{\rm if}~n \equiv 3(\bmod\,4).
\end{cases}$$
For $m\equiv 0(\bmod\,4)$ there is defined an isomorphism
$$L^{4k}(\Z) \xymatrix{\ar[r]^-{\di{\cong}}&} \Z;\quad (C,\phi)
\mapsto \,{\rm signature}(H^{2k}(C),\phi_0).$$
\indent
A $CW$ structure on an oriented $n$--dimensional manifold with $\pi_1(N)=\rho$
and universal cover $\wwtilde{N}$ and the Alexander--Whitney--Steenrod
diagonal construction on the cellular complex $C(\wwtilde{N})$
determine an $n$--dimensional symmetric Poincar\'e complex
$(C(\wwtilde{N}),\phi)$ over $\Z[\rho]$ with
$$\phi_0 = [N]\cap - \co C(\wwtilde{N})^{n-*} \to C(\wwtilde{N}).$$
The Mishchenko {\it symmetric signature} of $N$ is the cobordism class
$$\sigma^*(N) = (C,\phi) \in L^n(\Z[\rho]).$$
For $n=4k$ the image of $\sigma^*(N)$ in $L^{4k}(\Z)=\Z$ is just the
usual signature of $N$.

For any rings with involution $\Lambda,\Lambda'$ there are
defined products
$$\begin{array}{l}
L^{n}(\Lambda) \otimes L^{n'}(\Lambda')
\to L^{n+n'}(\Lambda\otimes \Lambda');
\ (C,\phi) \otimes (C',\phi') \mapsto  (C\otimes C',\phi\otimes \phi'),\\[1ex]
L_{n}(\Lambda) \otimes L^{n'}(\Lambda')
\to L_{n+n'}(\Lambda\otimes \Lambda');
\ (C,\psi) \otimes (C',\phi') \mapsto  (C\otimes C',\psi\otimes \phi')
\end{array}$$
as in Ranicki \cite{ranicki3}. The tensor product of group rings is given by
$$\Z[\pi]\otimes \Z[\pi'] = \Z[\pi\times \pi'].$$

\begin{thm}[Symmetric $L$--theory surgery product formula,
Ranicki \cite{ranicki-one}]
\label{surgeryproduct-one}\qua\newline
{\rm (i)}\qua The symmetric signature of a product $N \times N'$
of an $n$--dimensional manifold $N$ and an
$n'$--dimensional manifold $N'$ is the product of the symmetric signatures
$$\sigma^*(N \times N') = \sigma^*(N)\otimes \sigma^*(N')
\in L^{n+n'}(\Z[\pi_1(N) \times \pi_1(N')]).$$
{\rm (ii)}\qua The product of an $m$--dimensional
normal map $(f,b)\co M \to X$ and an $n$--di\-men\-sion\-al manifold $N$
is an $(m+n)$--dimensional normal map
$$(g,c) = (f,b) \times 1\co M \times N \to X \times N$$
with surgery obstruction
$$\sigma_*(g,c) = \sigma_*(f,b)\otimes \sigma^*(N) \in
L_{m+n}(\Z[\pi_1(X) \times \pi_1(N)]).$$
\end{thm}
\begin{proof} These formulae already hold on the chain homotopy
level, and chain equivalent symmetric/quadratic Poincar\'e complexes
are cobordant. In somewhat greater detail:

(i)\qua The symmetric Poincar\'e complex of a product $N''=N \times N'$ is
the product of the symmetric Poincar\'e complexes of $N$ and $N'$
$$(C(\wwtilde{N}''),\phi'') = 
(C(\wwtilde{N})\otimes C(\wwtilde{N}'),\phi\otimes\phi').$$
(ii)\qua The kernel quadratic Poincar\'e complex
$({\mathcal C}(g^!),\psi_c)$ of the product
normal map $(g,c)=(f,b)\times 1\co M \times N \to X \times N$
is the product of the kernel quadratic Poincar\'e complex
$({\mathcal C}(f^!),\psi_b)$
of $(f,b)$ and the symmetric Poincar\'e complex $(C(\wwtilde{N}),\phi)$
of $N$
$$({\mathcal C}(g^!),\psi_c) = 
({\mathcal C}(f^!)\otimes C(\wwtilde{N}),\psi_b \otimes \phi).\proved$$
\end{proof}

\begin{thm}
{\rm (i)}\qua {\rm (Shaneson \cite{Sh}, Wall \cite{cW-seventy})}\qua
The quadratic $L$--groups of $\Z[\Z^n]$ are given by
$$L_m(\Z[\Z^n]) = \sum\limits^n_{r=0}{n \choose r}L_{m-r}(\Z)~~(m \geqslant 0),$$
interpreting $L_{m-r}(\Z)$ for $m-r<0$ as $L_{m-r+4*}(\Z)$.\\
{\rm (ii)}\qua {\rm (Milgram and Ranicki \cite{MR}, Ranicki
\cite[Section~19]{ranicki-two})}\qua
The symmetric $L$--groups of $\Z[\Z^n]$ are given by
$$L^m(\Z[\Z^n]) = \sum\limits^n_{r=0}{n\choose r}L^{m-r}(\Z)~~(m \geqslant 0)$$
interpreting $L^{m-r}(\Z)$ for $m<r$ as
$$L^{m-r}(\Z) = \begin{cases}
0&{\rm if}~m=r-1,r-2\\[1ex]
L_{m-r}(\Z)&{\rm if}~m<r-2.
\end{cases}$$
\end{thm}

The symmetric signature of $T^n$
$$\sigma^*(T^n) = (C(\wwtilde{T}^n),\phi) = (0,\dots,0,1)
\in L^n(\Z[\Z^n]) = \sum\limits^n_{r=0}{n \choose r} L^{n-r}(\Z),$$
is the cobordism class of the $n$--dimensional symmetric Poincar\'e complex
$(C(\wwtilde{T}^n),\phi)$ over $\Z[\Z^n]$ with
$$C(\wwtilde{T}^n) = \bigotimes\limits_n C(\wtilde{S}^1),~
{\rm rank}_{\Z[\Z^n]}C_r(\wwtilde{T}^n) = {n \choose r}.$$
The surgery obstruction
$$E_8 \times T^n  = (0,\dots,0,1) \in
L_n(\Z[\Z^n]) = \sum\limits^n_{r=0}{n \choose r}L_{n-r}(\Z)$$
is the cobordism class of the $n$--dimensional quadratic Poincar\'e complex
over $\Z[\Z^n]$
$$(C,\psi) = (\Z^8,E_8)\otimes (C(\wwtilde{T}^n),\phi)$$
with
$${\rm rank}_{\Z[\Z^n]}C_r = 8 {n \choose r}.$$

\section{Almost $(-1)^n$--symmetric forms}\label{S-almost}

The surgery obstruction of the $(4m+2n)$--dimensional normal map
$$(f,b) = (f_0,b_0)\times 1\co M_0^{4m} \times T^{2n} \to S^{4m} \times T^{2n}$$
is given by the instant surgery obstruction of \fullref{S-instant}
and the surgery product formula of \fullref{S-product} to be the
Witt class
$$\sigma_*(f,b) = (K,\lambda,\mu) \in L_{4m+2n}(\Z[\Z^{2n}])$$
of the instant form $(K,\lambda,\mu)$ of the $2n$--dimensional
quadratic Poincar\'e complex
$$(C,\psi) = (\Z^8,E_8) \otimes (C(\wwtilde{T}^{2n}),\phi),$$
with
$${\rm rank}_{\Z[\Z^{2n}]}K = 
8\, {\rm rank}_{\Z[\Z^{2n}]}C_n(\wwtilde{T}^{2n}) =  8{2n
\choose n}.$$
In principle, it is possible to compute $(\lambda,\mu)$ directly from
the $(4m+2n)$--dimensional symmetric Poincar\'e complex $E_8 \otimes
(C(\wwtilde{T}^n),\phi)$.  In practice, we shall use the almost
symmetric form surgery product formula of Clauwens
\cite{clauwens-one,clauwens2,clauwens3}, which is the
analogue for symmetric Poincar\'e complexes of the instant surgery
obstruction of \fullref{S-instant}.  We establish a product formula for
almost symmetric forms which will be used in \fullref{S-almosttorus} to
obtain an almost $(-1)^n$--symmetric form for $T^{2n}$ of rank
$2^n\leqslant {2n \choose n}$, and hence a representative
$(-1)^n$--quadratic form for $\sigma_*(f,b)\in L_{4m+2n}(\Z[\Z^{2n}])$ of
rank $2^{n+3} \leqslant 8{2n \choose n}$.

\begin{defn}
Let $R$ be a ring with involution.\\
(i)\qua An {\it almost $(-1)^n$--symmetric form} $(A,\alpha)$
over $R$ is a f.g.\ free $R$--module $A$ together with a nonsingular pairing
$\alpha\co A \to A^*$ such that the endomorphism
$$\beta = 1+(-1)^{n+1}\alpha^{-1}\alpha^*\co  A \to A$$
is nilpotent, i.e. $\beta^N=0$ for some $N \geqslant 1$.\\
(ii)\qua A {\it sublagrangian} of an almost $(-1)^n$--symmetric form
$(A,\alpha)$ is a direct summand $L \subset A$ such that $L\subseteq
L^{\perp}$, where
$$L^{\perp} := \{b \in A\,\vert\,\alpha(b)(A)=\alpha(A)(b)=\{0\}\}.$$
A {\it lagrangian} is a sublagrangian $L$ such that
$$L = L^{\perp}.$$
(iii)\qua The {\it almost $(-1)^n$--symmetric Witt group} $AL^{2n}(R)$
is the abelian group of isomorphism classes of almost $(-1)^n$--symmetric
forms $(A,\alpha)$ over $R$ with relations
$$(A,\alpha) = 0~\hbox{\rm if $(A,\alpha)$ admits a lagrangian}$$
and addition by
$$(A,\alpha) + (A',\alpha')  = (A \oplus A',\alpha \oplus \alpha').$$
\end{defn}

\begin{example}
A nonsingular $(-1)^n$--symmetric form $(A,\alpha)$ is an
almost $(-1)^n$--sym\-met\-ric form such that
$$\alpha = (-1)^n\alpha^*\co A \to A^*$$
so that $1+(-1)^{n+1}\alpha^{-1}\alpha^*=0\co A \to A$.
\end{example}

An almost $(-1)^n$--symmetric form $(R^q,\alpha)$ on a f.g.\  free
$R$--module of rank $q$ is represented by an invertible $q \times
q$ matrix $\alpha=(\alpha_{rs})$ such that the $q \times q$ matrix
$$1+(-1)^{n+1} \alpha^{-1}\alpha^*$$
is nilpotent.

\begin{defn} The {\it instant form} of a $2n$--dimensional
symmetric Poincar\'e complex $(C,\phi)$ over $R$
is the almost $(-1)^n$--symmetric form over $R$
$$(A,\alpha) = 
\bigg({\rm coker}\biggl(\begin{pmatrix}
d^* & 0 \\ -\phi^*_0 & d \end{pmatrix}\co 
C^{n-1} \oplus C_{n+2}  \to C^n \oplus C_{n+1}\biggr),
\bigg[\begin{matrix} \phi_0+d\phi_1 & d \\
(-1)^nd^* & 0 \end{matrix}\bigg]\bigg).$$
\end{defn}

\begin{example} \label{iso}
If $\phi_0\co C^{2n-*} \to C$ is an isomorphism
the instant almost $(-1)^n$--sym\-met\-ric form is
$$(A,\alpha) = (C^n,\phi_0+d\phi_1).$$
\end{example}

Every $2n$--dimensional symmetric Poincar\'e complex $(C,\phi)$
over a ring with involution $R$
is chain equivalent to a complex $(C',\phi')$ such that
 $\phi'_0\co {C'}^{2n-*} \to C'$ is an isomorphism, with
 $$C'\co C^0 \to C^1 \to \dots \to C^{n-1} \to A^* \to C_{n-1} \to
 \dots \to C_1 \to C_0$$
and
$$\phi'_0+d'\phi'_1 = \alpha\co {C'}^n = A \to  C'_n = A^*.$$
(We shall not actually need this chain equivalence, since
$\phi_0\co C^{2n-*} \to C$ is an isomorphism for $C=C(\wwtilde{T}^{2n})$,
so \fullref{iso} will apply). The instant form defines a forgetful map
$$L^{2n}(R) \to AL^{2n}(R);\quad (C,\phi) \mapsto (A,\alpha).$$

\begin{proposition}[Ranicki {{\cite[36.3]{ranicki-four}}}]
The almost $(-1)^n$--symmetric Witt group of $\Z$ is given by
$$AL^{2n}(\Z) = \begin{cases}
\Z&\hbox{\it if $n \equiv 0(\bmod\,2)$}\\
0&\hbox{\it if $n \equiv 1(\bmod\,2)$}
\end{cases}$$
with $L^{4k}(\Z) \to AL^{4k}(\Z)$ an isomorphism.
The Witt class of an almost symmetric form $(A,\alpha)$ over $\Z$ is
$$(A,\alpha) = {\rm signature}(\Q\otimes A,(\alpha+\alpha^*)/2)
\in AL^{4k}(\Z) = L^{4k}(\Z) = \Z.$$
\end{proposition}

The almost $(-1)^n$--symmetric $L$--group $AL^{2n}(R)$ was denoted
$LAsy^0_{h,S_{(-1)^n}}(R)$ in~\cite{ranicki-four}.

\begin{defn} \label{asym}
The {\it almost symmetric signature} of a $2n$--dimensional
manifold $N^{2n}$ with $\pi_1(N)=\rho$ is the Witt class
$$\sigma^*(N) = (A,\alpha) \in AL^{2n}(\Z[\rho])$$
of the instant almost $(-1)^n$--symmetric form $(A,\alpha)$ over $\Z[\rho]$
of the $2n$--dimensional symmetric Poincar\'e complex
$(C(\wwtilde{N}),\phi)$ over $\Z[\rho]$.
\end{defn}

The forgetful map $L^{2n}(\Z[\rho]) \to AL^{2n}(\Z[\rho])$ sends the
symmetric signature $\sigma^*(N) \in L^{2n}(\Z[\rho])$ to the almost
symmetric signature $\sigma^*(N) \in AL^{2n}(\Z[\rho])$.

For any rings with involution $R_1,R_2$ there is defined a product
$$\begin{array}{l}
AL^{2n_1}(R_1) \otimes AL^{2n_2}(R_2) \to AL^{2n_1+2n_2}(R_1 \otimes R_2);\\[1ex]
\hphantom{AL^{2n_1}(R_1) \otimes AL^{2n_2}(R_2) \to }
(A_1,\alpha_1) \otimes (A_2,\alpha_2) \mapsto (A_1\otimes A_2,\alpha_1
\otimes \alpha_2).
\end{array}$$

\begin{proposition}\label{surgeryproduct-three}
The almost symmetric signature of a product $N=N_1\times
N_2$ of $2n_i$--dimensional manifolds $N_i$ with
$\pi_1(N_i)=\rho_i$ and almost $(-1)^{n_i}$--symmetric forms
$(\Z[\rho_i]^{q_i},\alpha_i)$ $(i=1,2)$ is the product
$$\begin{array}{ll}
\sigma^*(N_1 \times N_2)&= \sigma^*(N_1)\otimes \sigma^*(N_2)\\[1ex]
&\in {\rm im}(AL^{2n_1}(\Z[\rho_1]) \otimes
AL^{2n_2}(\Z[\rho_2]) \to AL^{2n_1+2n_2}(\Z[\rho_1 \times \rho_2])).
\end{array}$$
\end{proposition}
\begin{proof}
The almost $(-1)^{n_1+n_2}$--symmetric form  $(A,\alpha)$
of $N_1 \times N_2$ is defined on
$$A = C^{n_1+n_2}(\wwtilde{N}_1 \times \wwtilde{N}_2)
 = \bigoplus\limits_{(p_1,p_2) \in S}
C^{p_1}(\wwtilde{N}_1)\otimes C^{p_2}(\wwtilde{N}_2)$$
with
$$S = \{(p_1,p_2) \,\vert\,p_1+p_2=n_1+n_2\}.$$
Define an involution
$$T\co S \to S;\quad (p_1,p_2) \mapsto (2n_1-p_1,2n_2-p_2),$$
and let $U \subset S\backslash \{(n_1,n_2)\}$ be any subset such that
$S$ decomposes as a disjoint union
$$S = \{(n_1,n_2)\} \cup U \cup T(U).$$
The submodule
$$L = \bigoplus\limits_{(p_1,p_2) \in U}
C^{p_1}(\wwtilde{N}_1)\otimes C^{p_2}(\wwtilde{N}_2)
\subseteq A$$
is a sublagrangian of $(A,\alpha)$ such that
$$(L^{\perp}/L,[\alpha]) = (C^{n_1}(\wwtilde{N}_1),\alpha_1)\otimes
(C^{n_2}(\wwtilde{N}_2),\alpha_2).$$
The submodule
$$\Delta_{L^{\perp}} = \{(b,[b])\,\vert\,b \in L^{\perp}\} \subset A
\oplus (L^{\perp}/L)$$
is a lagrangian of $(A,\alpha) \oplus (L^{\perp}/L,-[\alpha])$,  and
$$\begin{array}{ll}
(A,\alpha)&= (L^{\perp}/L,[\alpha]) = 
(C^{n_1}(\wwtilde{N}_1),\alpha_1)\otimes
(C^{n_2}(\wwtilde{N}_2),\alpha_2)\\[1ex]
&\in {\rm im}(AL^{2n_1}(\Z[\rho_1]) \otimes
AL^{2n_2}(\Z[\rho_2]) \to AL^{2(n_1+n_2)}(\Z[\rho_1 \times \rho_2])).
\end{array}$$
This completes the proof.
\end{proof}
The product of a nonsingular $(-1)^m$--quadratic form
$(K,\lambda,\mu)$ over $\Lambda$ and a $2n$--dimensional symmetric
Poincar\'e complex $(C,\phi)$ over $R$ is a $2(m+n)$--dimensional
quadratic Poincar\'e complex $(K_{*-m}\otimes C,(\lambda,\mu)\otimes \phi)$
over $\Lambda'=\Lambda\otimes R$, as in Ranicki \cite{ranicki-one},
with $K_{*-m}$ the $2m$--dimensional f.g.\ free $\Lambda$--module chain complex
concentrated in degree $m$
$$K_{*-m}\co 0 \to \dots \to 0 \to K \to 0 \to \dots \to 0.$$
The pairing
$$L_{2m}(\Lambda) \otimes L^{2n}(R) \to L_{2m+2n}(\Lambda\otimes R);
\ (K,\lambda,\mu) \otimes (C,\phi) \mapsto (K_{*-m}\otimes C,(\lambda,\mu)\otimes\phi)$$
has the following generalization.

\begin{defn}
The {\it product} of a nonsingular $(-1)^m$--quadratic form
$(K,\lambda,\mu)$ over $\Lambda$ and an almost $(-1)^n$--symmetric form
$(A,\alpha)$ over $R$ is  the nonsingular
$(-1)^{m+n}$--quadratic form over $\Lambda'=\Lambda\otimes R$
$$(K',\lambda',\mu') = (K,\lambda,\mu)\otimes (A,\alpha)$$
with
$$K' = K \otimes A,~(\lambda',\mu') = (\lambda,\mu)\otimes\alpha = 
(\psi'+(-1)^{m+n}{\psi'}^*,\psi')$$
determined by the $\Lambda'$--module morphism
$$\psi' = \psi \otimes \alpha\co 
K' = K \otimes A \to {K'}^* = K^*\otimes A^*$$
with $\psi\co K \to K^*$ a $\Lambda$--module morphism such that
$$(\lambda,\mu) = (\psi+(-1)^m\psi^*,\psi).$$
\end{defn}

In particular, if $K=\Lambda^p$ then $\psi$ is given by a
$p \times p$ matrix $\psi=\{\psi_{ij}\}$ over $\Lambda$, and if $A=R^q$
then $\alpha=\{\alpha_{rs}\}$ is given by a $q \times q$ matrix
over $R$, so that
$$\psi' = \psi \otimes \alpha$$
is the $pq\times pq$ matrix over $\Lambda'$ with
$$\psi'_{tu} = \psi_{ij}\otimes \alpha_{rs}~{\rm if}~t=(i-1)p+r,~u=(j-1)p+s.$$
\indent If $(A,\alpha)$ is an almost $(-1)^n$--symmetric form over $R$
with a sublagrangian $L \subset A$ the induced
almost $(-1)^n$--symmetric form $(L^{\perp}/L,[\alpha])$ over $R$ is
such that
$$\Delta_{L^{\perp}} = \{(b,[b])\,\vert\,b \in L^{\perp}\} \subset A
\oplus (L^{\perp}/L)$$
is a lagrangian of $(A,\alpha) \oplus (L^{\perp}/L,-[\alpha])$,  and
$$(K,\lambda,\mu)\otimes (A,\alpha) = 
(K,\lambda,\mu)\otimes (L^{\perp}/L,[\alpha]) \in L_{2m+2n}(\Lambda').$$
In particular, if $L$ is a lagrangian of $(A,\alpha)$ then
$$(K,\lambda,\mu)\otimes (A,\alpha) = 0 \in L_{2m+2n}(\Lambda'),$$
so that the product
$$L_{2m}(\Lambda) \otimes AL^{2n}(R) \to L_{2m+2n}(\Lambda');\quad 
(K,\lambda,\mu) \otimes (A,\alpha) \mapsto (K\otimes A,(\lambda,\mu)\otimes\alpha)$$
is well-defined.

\begin{thm}[Almost symmetric $L$--theory surgery product formula, Clauwens
\cite{clauwens-one}]
\label{surgeryproduct-two}
{\rm (i)}\qua The product
$$L_{2m}(\Lambda) \otimes L^{2n}(R) \to L_{2m+2n}(\Lambda\otimes R);
\ (K,\lambda,\mu) \otimes (C,\phi) \mapsto (K_{*-m}\otimes C,(\lambda,\mu)\otimes\phi)$$
factors through the product
$$L_{2m}(\Lambda) \otimes AL^{2n}(R) \to L_{2m+2n}(\Lambda\otimes R);
\ (K,\lambda,\mu) \otimes (A,\alpha) \mapsto (K\otimes A,(\lambda,\mu)\otimes\alpha).$$
{\rm (ii)}\qua Let $(f,b)\co M \to X$ be a $2m$--dimensional normal map with surgery obstruction
$$\sigma_*(f,b) = (\Z[\pi]^p,\lambda,\mu) \in L_{2m}(\Z[\pi])~~(\pi=\pi_1(X)),$$
and let $N$ be a $2n$--dimensional manifold with almost $(-1)^n$--symmetric
signature
$$\sigma^*(N) = (\Z[\rho]^q,\alpha) \in AL^{2n}(\Z[\rho])~~(\rho=\pi_1(N)).$$
The surgery obstruction of the $(2m+2n)$--dimensional normal map
$$(g,c) = (f,b) \times 1\co M \times N \to X \times N$$
is given by
$$\begin{array}{ll}
\sigma_*(g,c)&= (\Z[\pi \times \rho]^{pq},(\lambda,\mu)\otimes \alpha)\\[1ex]
&\in {\rm im}(L_{2m}(\Z[\pi]) \otimes
AL^{2n}(\Z[\rho]) \to L_{2m+2n}(\Z[\pi \times \rho])).
\end{array}$$
{\rm (iii)}\qua The surgery obstruction of the product
$2(m+n_1+n_2)$--dimensional normal map
$$(g,c) = (f,b) \times 1\co M \times N_1 \times N_2 \to X \times N_1
\times N_2$$
is given by
$$\sigma_*(g,c) = (\Z[\pi\times \rho_1 \times \rho_2]^{pq_1q_2},
(\lambda,\mu)\otimes \alpha_1 \otimes \alpha_2) \in
L_{2(m+n_1+n_2)}(\Z[\pi\times \rho_1 \times \rho_2]).$$
\end{thm}
\begin{proof}
(i)\qua By construction.\\
(ii)\qua It may be assumed that $(f,b)\co M \to X$ is an $m$--connected
$2m$--dimensional normal map, with kernel $(-1)^m$--quadratic
form over $\Z[\pi]$
$$(K_m(M),\lambda,\mu) = (\Z[\pi]^p,\lambda,\mu).$$
The product $(g,c)=(f,b)\times 1\co M \times N \to X \times N$ is $m$--connected,
with quadratic Poincar\'e complex
$$(C,\psi) = (K_m(M),\lambda,\mu) \otimes (C(\wwtilde{N}),\phi)$$
and kernel $\Z[\pi \times \rho]$--modules
$$K_*(M \times N) = K_m(M)\otimes H_{*-m}(\wwtilde{N}).$$
Let $(f',b')\co M' \to X \times N$ be the bordant $(m+n)$--connected
normal map obtained from $(g,c)$ by surgery below the middle
dimension, using $(C,\psi)$ as in \fullref{below} (ii).
The kernel $(-1)^{m+n}$--quadratic form over $\Z[\pi\times \rho]$ of $(f',b')$
is the instant form of $(C,\psi)$, which
is just the product of $(K_m(M),\lambda,\mu)$ and
the almost $(-1)^n$--symmetric form $(\Z[\rho]^q,\alpha)$
\begin{multline*}
(K_{m+n}(M'),\lambda',\mu') =\\[1ex]
\bigg({\rm coker}\biggl(\begin{pmatrix}
d^* & 0 \\ (-1)^{m+n+1}(1{+}T)\psi_0 & d \end{pmatrix}\co 
C^{m+n-1} \oplus C_{m+n+2}  \to C^{m+n} \oplus C_{m+n+1}\biggr),\\[1ex]
\hskip200pt \bigg[\begin{matrix} \psi_0+(-1)^{m+n}\psi^*_0 & d \\
(-1)^{m+n}d^* & 0 \end{matrix}\bigg],\bigg[\begin{matrix} \psi_0 & d \\
0 & 0 \end{matrix}\bigg]\bigg)\\[2ex]
\hphantom{(K_{m+n}(M'),\lambda',\mu')~}
= (\Z[\pi \times \rho]^{pq},(\lambda,\mu)\otimes \alpha).
\end{multline*}
The surgery obstruction of $(g,c)$ is thus given by
$$\begin{array}{ll}
\sigma_*(g,c)&= \sigma_*(f',b') = (K_{m+n}(M'),\lambda',\mu')\\[1ex]
&= (\Z[\pi \times \rho]^{pq},(\lambda,\mu)\otimes \alpha)
\in L_{2m+2n}(\Z[\pi \times \rho]).
\end{array}$$
(iii)\qua Combine (i) and (ii) with \fullref{surgeryproduct-three}.
\end{proof}

\section{The almost $(-1)^n$--symmetric form of $\mathbf{T^{2n}}$}
\label{S-almosttorus}

Geometrically, $--\times T^{2n}$ sends the surgery obstruction
$\sigma_*(f_0,b_0)=E_8 \in L_{4m}(\Z)$ to the surgery obstruction
$$E_8\times T^{2n}~ = 
\sigma_*(f_n,b_n) \in L_{4m+2n}(\Z[\Z^{2n}])$$
of the $(4m+2n)$--dimensional normal map
$$(f_n,b_n) = (f_0,b_0)\times 1\co M_0^{4m}\times T^{2n} \to
S^{4m}\times T^{2n}$$
given by product with the almost symmetric signature of
$$\begin{array}{ll}
T^{2n}&= S^1 \times S^1 \times \dots \times S^1~~(2n~{\rm factors})\\[1ex]
&= T^2 \times T^2 \times \dots \times T^2~~(n~{\rm factors}).
\end{array}$$
In order to apply the almost symmetric surgery product formula
(see \fullref{surgeryproduct-two}) for
$N^{2n}=T^{2n}$ it therefore suffices to work out the almost
$(-1)$--symmetric form $(C^1(\wwtilde{T}^2),\alpha)$ of $T^2$.

The symmetric Poincar\'e structure $\phi=\{\phi_s\vert s \geqslant 0\}$
of the universal cover $\wtilde{S}^1=\R$ of $S^1$ is given by
$$\begin{array}{l}
d = 1-z\co C_1(\R) = \Z[z,z^{-1}] \to C_0(\R) = \Z[z,z^{-1}],\\[1ex]
\phi_0 = \begin{cases}
1\co  C^0(\R) = \Z[z,z^{-1}] \to C_1(\R) = \Z[z,z^{-1}]\\[0.5ex]
z\co  C^1(\R) = \Z[z,z^{-1}] \to C_0(\R) = \Z[z,z^{-1}],
\end{cases}\\[3ex]
\phi_1 = -1\co C^1(\R) = \Z[z,z^{-1}] \to C_1(\R) = \Z[z,z^{-1}].
\end{array}$$
Write
$$\Lambda = \Z[\pi_1(T^2)] = \Z[z_1,z_1^{-1},z_2,z_2^{-1}].$$
The Poincar\'e duality of $\wwtilde{T}^2=\R^2$ is the $\Lambda$--module
chain isomorphism given by the chain-level K\"unneth formula to be
$$\xymatrix@C+50pt@R+40pt
{C(\wwtilde{T}^2)^{2-*} \ar@<-5ex>[d]_-{\displaystyle{\phi_0}}\co 
C^0{=}\Lambda \ar[r]^-{\displaystyle{d^*{=}\begin{pmatrix}
z_2{-}1 \\ 1{-}z_1^{-1} \end{pmatrix}}} \ar@<5ex>[d]^{\displaystyle{1}} &
C^1{=}\Lambda{\oplus}\Lambda\ar[r]^-{\displaystyle{d^*{=}(1{-}z^{-1}_1~1{-}z_2)}}
\ar[d]^-{\displaystyle{\begin{pmatrix} 0 & -z_1 \\ z^{-1}_2 & 0 \end{pmatrix}}} &
C^2{=}\Lambda \ar[d]^-{\displaystyle{-z_1z^{-1}_2}} \\
C(\wwtilde{T}^2)\co  C_2{=}\Lambda
\ar[r]^-{\displaystyle{d{=}\begin{pmatrix} 1{-}z_1 \\ 1{-}z_2^{-1}
\end{pmatrix}}} & C_1{=}\Lambda{\oplus}\Lambda
\ar[r]^-{\displaystyle{d{=}(z_2^{-1}{-}1~1{-}z_1)}} & C_0 {=} \Lambda.}
$$
The chain homotopy
$$\phi_1\co \phi_0\simeq T\phi_0\co C(\wwtilde{T}^2)^{2-*} \to
C(\wwtilde{T}^2)$$
is given by
$$\phi_1 = 
\begin{cases}
\begin{pmatrix} 1 & -z_2 \end{pmatrix}\co C^1 = \Lambda\oplus \Lambda \to C_2 = \Lambda\\[1ex]
\begin{pmatrix} -z_1 \\ 1 \end{pmatrix}\co C^2 = \Lambda \to C_1 = \Lambda\oplus \Lambda.
\end{cases}$$

\begin{proposition}\label{t-two}
The almost $(-1)$--symmetric form of $T^2$ is given by
$(C^1,\alpha)$ with
$$\alpha = \phi_0-\phi_1d^* = 
\begin{pmatrix} 1-z_1 & z_1z_2-z_1-z_2\\
1 & 1-z_2 \end{pmatrix}\co 
C^1 = \Lambda \oplus \Lambda \to C_1 = \Lambda \oplus \Lambda.$$
\end{proposition}
\begin{proof} By construction, noting that
$$\begin{array}{ll}
1+\alpha^{-1}\alpha^*&= 
\begin{pmatrix} -(1-z_1)(1-z_2^{-1}) & z_1(1-z_2)(1-z_2^{-1}) \\[0.5ex]
-z_2^{-1}(1-z_1)(1-z_1^{-1}) & (1-z_1)(1-z_2^{-1}) \end{pmatrix}\co \\[2ex]
&\hskip100pt
C^1 = \Lambda \oplus \Lambda \to C^1 = \Lambda \oplus \Lambda
\end{array}$$
is nilpotent, with
$$(1+\alpha^{-1}\alpha^*)^2 = 0\co 
C^1 = \Lambda \oplus \Lambda \to C^1 = \Lambda \oplus \Lambda.\proved$$
\end{proof}

\begin{rem}
An almost $(-1)^n$--symmetric form $(R^q,\alpha)$
over $R$ determines a nonsingular $(-1)^n$--quadratic form
$(R[1/2]^q,\lambda,\mu)$ over $R[1/2]$, with
$$\begin{array}{l}
\lambda(x,y) = (\alpha(x,y)+(-1)^n\overline{\alpha(y,x)})/2,~
\mu(x) = \alpha(x)(x)/2.
\end{array}$$
In particular, the almost $(-1)$--symmetric form $(\Lambda\oplus
\Lambda,\alpha)$ of $T^2$ determines the nonsingular $(-1)$--quadratic
form $(\Lambda[1/2]\oplus \Lambda[1/2],\lambda,\mu)$ over
$\Lambda[1/2]=\Z[\Z^2][1/2]$, with
$$\begin{array}{ll}
\lambda&= (\alpha-\alpha^*)/2\\[1ex]
&= \begin{pmatrix} ((z_1)^{-1}-z_1)/2 & (1-z_1z_2-z_1-z_2)/2 \\[1ex]
(-1+(z_1)^{-1}(z_2)^{-1}+(z_1)^{-1}+(z_2)^{-1})/2 & ((z_2)^{-1}-z_2)/2 \end{pmatrix}
\end{array}$$
the invertible skew-symmetric $2 \times 2$ matrix exhibited in
\cite[Example, p120]{G}.
\end{rem}

\section{An explicit form representing $\mathbf{E_8\times T^{2n} \in L_{4*+2n}(\Z[\Z^{2n}])}$}\label{S-explicit}

Write the generators of the free abelian group
$\pi_1(T^{2n})=\Z^{2n}$ as
$$z_1,z_2,\ldots,z_{2n-1},z_{2n},$$
so that
$$\Z[\Z^{2n}] = \Z[z_1,z_1^{-1},z_2,z_2^{-1},\dots,z_{2n},z_{2n}^{-1}].$$
The expression of $T^{2n}$ as an $n$--fold cartesian product of
$T^2$'s
$$T^{2n} = T^2 \times T^2 \times \dots \times T^2$$
gives
$$\Z[\Z^{2n}] = \Z[z_1,z_1^{-1},z_2,z_2^{-1}]\otimes
\Z[z_3,z_3^{-1},z_4,z_4^{-1}]\otimes\dots
\otimes\Z[z_{2n-1},z_{2n-1}^{-1},z_{2n},z_{2n}^{-1}].$$
For $i=1,2,\dots,n$ define the invertible $2 \times 2$ matrix over
$\Z[z_{2i-1},z_{2i-1}^{-1},z_{2i},z_{2i}^{-1}]$
$$\alpha_i = \begin{pmatrix} 1-z_{2i-1} & z_{2i-1}z_{2i}-z_{2i-1}-z_{2i} \\
1 & 1-z_{2i} \end{pmatrix}.$$
The generator $1=E_8 \in L_0(\Z)=\Z$ is represented by the nonsingular
quadratic form $(\Z^8,\psi_0)$ over $\Z$ with:
$$\psi_0 = \begin{pmatrix}
1&0&0&1&0&0&0&0 \\
0&1&1&0&0&0&0&0 \\
0&0&1&1&0&0&0&0 \\
0&0&0&1&1&0&0&0 \\
0&0&0&0&1&1&0&0 \\
0&0&0&0&0&1&1&0 \\
0&0&0&0&0&0&1&1 \\
0&0&0&0&0&0&0&1
\end{pmatrix}$$

\begin{thm} \label{T-E-eight}
The $2^{n+3} \times 2^{n+3}$ matrix over $\Z[\Z^{2n}]$
\[
\psi_n = \psi_0 \otimes \alpha_1 \otimes \alpha_2 \dots \otimes \alpha_n
\]
is such that
\[
E_8 \times T^{2n} = (\Z[\Z^{2n}]^{2^{n+3}},\psi_n) \in
L_{2n}(\Z[\Z^{2n}])\, .
\]
\end{thm}
\begin{proof} A direct application of the almost symmetric surgery
product formula (\fullref{surgeryproduct-two}), noting that $\alpha_1$,
$\alpha_2$, $\dots$, $\alpha_n$ are copies of the almost
$(-1)$--symmetric form of $T^2$ obtained in \fullref{t-two}.
\end{proof}

\section{Transfer invariance} \label{S-transfer}

A covering map $p\co T^n \to T^n$ induces an injection of the
fundamental group in itself
$$p_*\co \pi_1(T^n) = \Z^n \to \pi_1(T^n) = \Z^n$$
as a subgroup of finite index, say $q=[\Z^n\co p_*(\Z^n)]$.
Given a $\Z[\Z^n]$--module $K$ let $p^!K$ be the $\Z[\Z^n]$--module
defined by the additive group of $K$ with
$$\Z[\Z^n] \times p^!K \to p^!K;\quad (a,b) \mapsto p_*(a)b.$$
In particular
$$p^!\Z[\Z^n] = \Z[\Z^n]^q.$$
The restriction functor
$$p^!\co \{\Z[\Z^n]\hbox{\rm -modules}\} \to
\{\Z[\Z^n]\hbox{\rm -modules}\};\quad K \mapsto p^!K$$
induces transfer maps in the quadratic $L$--groups
$$p^!\co L_m(\Z[\Z^n]) \to L_m(\Z[\Z^n]);\quad (C,\psi) \mapsto p^!(C,\psi).$$
\begin{proposition} \label{invariance}
The image of the (split) injection
$$L_0(\Z) \to L_n(\Z[\Z^n])=\sum\limits^n_{r=0}{n \choose r}L_{n-r}(\Z);
\quad E_8 \mapsto E_8 \times T^n$$
is the subgroup of the transfer-invariant elements
$$L_n(\Z[\Z^n])^\mathit{INV} = \{x \in L_n(\Z[\Z^n])\,\vert\,
p^!x=x~\hbox{\it for all}~p\co T^n\to T^n\}.$$
\end{proposition}
\begin{proof} See Ranicki \cite[Chapter 18]{ranicki-two}.
\end{proof}

\begin{example}
(i)\qua Write
$$\Lambda = \Z[\Z^2] = \Z[z_1,z^{-1}_1,z_2,z^{-1}_2].$$
Here is an explicit verification that
$$p^!(E_8\times T^2) = E_8 \times T^2 \in L_2(\Lambda)$$
for the double cover
$$p\co T^2 = S^1   \times S^1 \to T^2;\quad (w_1,w_2) \mapsto ((w_1)^2,w_2)$$
with
$$p_*\co \pi_1(T^2) = \Z^2 \to \Z^2;\quad z_1 \mapsto (z_1)^2,~z_2 \mapsto z_2$$
the inclusion of a subgroup of index 2. For any $j_1,j_2 \in \Z$ the transfer of the
$\Lambda$--module morphism $z_1^{j_1}z_2^{j_2}\co \Lambda \to \Lambda$ is given by the
$\Lambda$--module morphism
$$p^! (z_1^{j_1} z_2^{j_2}) =
\begin{cases}
\begin{pmatrix} (z_1)^{j_1/2}z_2^{j_2} & 0 \\ 0 & (z_1)^{j_1/2}z_2^{j_2}
\end{pmatrix}\co &\\
\hskip25pt p^!\Lambda = \Lambda\oplus \Lambda \to p^!\Lambda = \Lambda\oplus \Lambda&
\hbox{\rm if $j_1$ is even}\\[2ex]
\begin{pmatrix} 0 & (z_1)^{(j_1+1)/2}z_2^{j_2} \\
(z_1)^{(j_1-1)/2}z_2^{j_2} & 0 \end{pmatrix}\co &\\
\hskip25pt p^!\Lambda = \Lambda\oplus \Lambda \to p^!\Lambda = \Lambda\oplus \Lambda&
\hbox{\rm if $j_1$ is odd}.
\end{cases}
$$
\indent The transfer of the almost $(-1)$--symmetric form of $T^2$ over $\Lambda$
$$(C^1(\wwtilde{T}^2),\alpha) = 
\Big(\Lambda \oplus \Lambda,\begin{pmatrix} 1-z_1 & z_1z_2-z_1-z_2 \\ 1 & 1-z_2 \end{pmatrix}\Big)$$
is the almost $(-1)$--symmetric form over $\Lambda$
$$p^!(C^1(\wwtilde{T}^2),\alpha) = \Big(\Lambda \oplus \Lambda \oplus \Lambda \oplus \Lambda,
\begin{pmatrix}
1 & -z_1 & -z_2 & z_1z_2-z_1 \\
-1 & 1 & z_2-1 & -z_2 \\
1 & 0 & 1-z_2 & 0 \\
0 & 1 & 0 & 1-z_2
\end{pmatrix}\Big)$$
The $\Lambda$--module morphisms
$$\begin{array}{l}
i = \begin{pmatrix} z_1-z_1z_2 \\ 0 \\ -z_1 \\ 1 \end{pmatrix}\co 
\Lambda \to \Lambda \oplus \Lambda \oplus \Lambda \oplus \Lambda,\\[5ex]
j = \begin{pmatrix} 1 & 0 & z_1-z_1z_2 \\
z_1^{-1} & 0 & 0 \\
0 & 1 & -z_1 \\
0 & 0 & 1 \end{pmatrix}\co 
\Lambda\oplus \Lambda \oplus \Lambda \to \Lambda \oplus \Lambda \oplus \Lambda \oplus \Lambda
\end{array}$$
are such that $i=j\vert_{0 \oplus 0 \oplus \Lambda}$ and
there is defined a (split) exact sequence
$$\xymatrix@C+10pt{0 \ar[r] & \Lambda\oplus \Lambda \oplus \Lambda \ar[r]^-{\displaystyle{j}} &
\Lambda \oplus \Lambda \oplus \Lambda \oplus \Lambda \ar[r]^-{\displaystyle{i^*p^!\alpha}} &
\Lambda \ar[r] & 0}$$
with
$$j^*(p^!\alpha)j = \begin{pmatrix} 1-z_1 & z_1z_2-z_1-z_2 & 0 \\
1 & 1-z_2 & 0 \\
0 & 0 & 0 \end{pmatrix}\co \Lambda \oplus \Lambda \oplus \Lambda \to \Lambda \oplus \Lambda \oplus \Lambda.$$
The submodule
$$L = i(\Lambda) \subset p^!(\Lambda \oplus \Lambda) = \Lambda \oplus \Lambda \oplus \Lambda \oplus \Lambda$$
is thus a sublagrangian of the almost $(-1)$--symmetric form
$p^!(C^1(\wwtilde{T}^2),\alpha)$ over $\Z[\Z^2]$ such that
$$(L^{\perp}/L,[p^!\alpha]) = (C^1(\wwtilde{T}^2),\alpha)$$
and
$$\begin{array}{ll}
p^!(E_8\times T^2)&= E_8 \otimes p^!(C^1(\wwtilde{T}^2),\alpha)\\[1ex]
&= E_8 \otimes (L^{\perp}/L,[p^!\alpha])\\[1ex]
&= E_8\otimes (C^1(\wwtilde{T}^2),\alpha) = E_8 \times T^2~ \in L_2(\Lambda).
\end{array}$$
(ii)\qua For any $n \geqslant 1$ replace $p$ by
$$p_n = p \times 1\co T^{2n} = T^2 \times T^{2n-2} \to T^{2n} = T^2 \times T^{2n-2}$$
to likewise obtain an explicit verification that
$$p_n^!(E_8\times T^{2n}) = E_8\times  T^{2n} \in L_{2n}(\Z[\Z^{2n}]).$$
\end{example}

\section{Controlled surgery groups}\label{S-control}

A {\em geometric $\Z[\pi]$--module\/} over a metric space $B$ is a
pair $(K, \varphi)$, where $K = \Z[\pi]^r$ is a free $\Z[\pi]$--module
with basis $S = \{e_1, \ldots, e_r \}$ and $\varphi \colon S \to B$
is a map.   The {\em $(\epsilon, \delta)$--controlled surgery group
$L_n (B; \Z, \epsilon, \delta)$}  (with trivial local fundamental
group) is defined as the group of $n$--dimensional
quadratic $\Z$--Poincar\'e complexes (see \cite{ranicki-one}) over $B$
of radius $< \delta$, modulo $(n+1)$--dimensional
quadratic $\Z$--Poincar\'e bordisms of radius $< \epsilon$.  Elements of
$L_{2n} (B; \Z, \epsilon, \delta)$ are represented by non-singular
$(-1)^n$--quadratic forms $(K, \lambda, \mu)$, where $K = \Z^r$ is
a geometric $\Z$--module over $B$, and $\lambda$
has radius $< \delta$, i.e., $\lambda (e_i, e_j) = 0$ if
$d(\varphi(e_i), \varphi(e_j)) \geqslant \delta$.   In matrix
representation $(K, \psi)$, this is equivalent to $\psi_{ij} = 0$ if
$d(\varphi(e_i), \varphi(e_j)) \geqslant \delta$. The radius of a bordism
is defined similarly.

In effect, Yamasaki \cite{mY87} defined an assembly map $H_n (B; \LL) \to
L_n (B; \Z, \epsilon, \delta)$, where $H_\ast (B; \LL)$ denotes
homology with coefficients in the 4--periodic simply-connected
surgery spectrum $\LL$ of Ranicki \cite[Chapter~25]{ranicki3}.

The following Stability Theorem is a key
ingredient in the construction of exotic ENR homology manifolds.

\begin{thm}[Stability; Pedersen, Quinn and Ranicki
\cite{PQR-zero-one}; Ferry \cite{sF-zero-three}; Pedersen and Yamasaki
\cite{PY}] \label{stability}\qua\\
Let $n\geqslant 0$ and suppose $B$ is a compact metric ENR.  Then there
exist constants $\epsilon_0>0$ and $\kappa>1$, which depend on $n$ and
$B$, such that the assembly map $H_n (B; \LL) \to L_n (B; \Z, \epsilon,
\delta)$ is an isomorphism if $\epsilon_0\geqslant \epsilon\geqslant
\kappa\delta$, so that
$$\varprojlim_{\epsilon}\varprojlim_{\delta}
L_n(B;\Z,\epsilon,\delta) = H_n (B;\mathbb{L}).$$
\end{thm}

We are interested in controlled surgery over the torus $T^{2n} =
\R^{2n} / \Z^{2n}$ equipped with the usual geodesic metric. Let
$(K,\psi)$ represent an element of $L_{2n} (\Z[\Z^{2n}])$, where
$K = \Z[\Z^{2n}]^r$.  Our next goal is to show that passing to a
sufficiently large covering space $p \colon T^{2n} \to T^{2n}$,
$(K, \psi)$ defines an element of $L_{2n} (T^{2n}; \Z, \epsilon,
\delta)$.  For simplicity, we assume that
\[
p_\ast \colon \pi_1 (T^{2n}) \cong \Z^{2n} \to \pi_1 (T^{2n})
\cong \Z^{2n}
\]
is given
by multiplication by $k > 0$, so that $p$ is a $k^{2n}$--sheeted
covering space.

Let $(\wwbar K, \wbar \psi) = \Z [\Z_k^{2n}] \otimes_{\Z[\Z^{2n}]} (K, \psi)$,
where the (right) $\Z[\Z^{2n}]$--module structure on $\Z[\Z_k^{2n}]$ is
induced by reduction modulo $k$.  The $\Z$--module $\wwtilde K$
underlying  $\wwbar K$ has basis $\Z_k^{2n} \times S$; if $g \in
\Z_k^{2n}$ and $e_i \in S$, we write $(g,e_i) = g \, e_i$.
Pick a point $x_0$ in the covering torus $T^{2n}$
viewed as a $\Z_k^{2n}$--space under the action of the group
of deck transformations.   Let
$\varphi (e_i) = x_0$, for every $e_i \in S$, and extend it
$\Z_k^{2n}$--equivariantly to obtain $\varphi \colon
\Z_k^{2n} \times S \to T^{2n}$.  Then, the pair $(\wwtilde K, \varphi)$
is a geometric $\Z$--module over $T^{2n}$ of dimension
$r k^{2n}$.

We now describe the quadratic $\Z$--module $(\wwtilde K, \wtilde \psi)$
induced by $(K, \psi)$ and the covering $p$.   Write
$$\wbar{\psi} =\sum_{g \in \Z_k^{2n}}  g \, \wbar{\psi}_g,$$
where each $\wbar{\psi}_g$ is a matrix with integer entries.
For basis elements $g e_i, f e_j \in \Z_k^{2n} \times S$, let
$\wtilde \psi (g e_i, f e_j) = \wbar{\psi}_{g^{-1} f} (e_i, e_j)$;  this
defines a bilinear $\Z$--form on the geometric $\Z$--module $\wwtilde K$.  For
a given quadratic $\Z[\Z^{2n}]$--module $(K, \psi)$, we show that
$(\wwtilde K, \wtilde \psi)$ has
diameter $< \delta$ over the (covering) torus $T^{2n}$, if $k$ is
sufficiently large.

Elements of $\Z^{2n}$  can be expressed uniquely as monomials
$$z^i = z_1^{i_1}\ldots z_{2n}^{i_{2n}}$$
where $i = (i_1, \ldots, i_{2n}) \in \Z^{2n}$ is
a multi-index.  We use the notation
$$| i | = \text{max} \, \{ | i_1 |, \ldots, | i_{2n} | \}.$$
  Any $z \in \Z[\Z^{2n}]$ can be expressed uniquely as
$$z = \sum_{i \in \Z^{2n}} \alpha_i \, z^i,$$
where $\alpha_i \in \Z$ is zero for all but finitely many values of $i$.
We define the {\em order of $z$} to be
$$o(z) = \text{max}\, \{ |i| \, : \, \alpha_i \ne 0 \}$$
and let
$$|\psi| = \text{max}\, \{o(\psi_{ij}), 1 \leqslant i, j \leqslant r \}.$$
Then,
$(\wwtilde K, \wtilde \psi)$ is a quadratic $\Z$--module over $T^{2n}$ of
radius $< \delta$, provided that $k > 2 | \psi | / \delta$.  Similarly,
quadratic $\Z[\Z^{2n}]$--Poincar\'e bordisms
induce quadratic $\Z$--Poincar\'e $\epsilon$--bordisms for $k$ large.

\subsection{The forgetful map}  We give an algebraic description of
the forget-control map
\[
\forget  \colon L_{2n}(T^{2n}; \Z, \epsilon, \delta) \to
L_{2n}(\Z[\Z^{2n}]),
\]
for $\epsilon$ and $\delta$ small.  Let $\sigma \in L(T^{2n}; \Z,
\epsilon, \delta)$ be represented by the $(-1)^n$--quadratic $\Z$--module
$(K, \psi)$ over $T^{2n}$ of radius $< \delta$, where $K$ has basis
$S= \{ e_1, \ldots, e_r \}$ and projection $\varphi \colon S \to T^{2n}$.
Consider the free $\Z[\Z^{2n}]$--module $\wwtilde K$ of rank $r$ generated by
$\wtilde S = \{ \wtilde{e}_1, \ldots, \wtilde{e}_r \}$ and let $\wtilde \varphi
\colon \wtilde S \to \R^{2n}$ be a map satisfying $q \circ \wtilde
\varphi (\wtilde{e}_i) = \varphi (e_i)$, $1 \leqslant i \leqslant r$, where
$q \colon \R^{2n} \to T^{2n} =  \R^{2n} / \Z^{2n}$ is the universal
cover.  If $\psi_{ij} \ne 0$ and
$\delta$ is small, there is a unique element $g_{ij}$ of $\Z^{2n}$
such that $d \left( \wtilde{\varphi} (\wtilde{e}_j) + g_{ij}, \wtilde{\varphi}
(\wtilde{e}_i) \right) < \delta$, where $d$ denotes Euclidean distance.
Let $\wtilde \psi = \left( \wtilde{\psi}_{ij} \right)$, ${1 \leqslant i, j \leqslant r}$ be the matrix
whose entries in $\Z[\Z^{2n}]$ are
\begin{equation}
\wtilde{\psi}_{ij} =
\begin{cases}
0, & \text{if $\psi_{ij} = 0 \,$;}  \\
\psi_{ij} \, g_{ij}, & \text{if $\psi_{ij} \ne 0\,$.}
\end{cases}
\end{equation}
The quadratic $\Z[\Z^{2n}]$--module $(\wwtilde K, \wtilde \psi)$ represents
$\forget (\sigma) \in L_{2n} (\Z[\Z^{2n})$.  Likewise,
quadratic $\Z$--Poincar\'e $\epsilon$--bordisms over $T^{2n}$
induce quadratic $\Z[\Z^{2n}]$--Poincar\'e bordisms.

\subsection{Controlled $\mathbf{E_8}$ over $\mathbf{T^{2n}}$}
\label{S-E-eight}
Starting with the $(-1)^n$--quadratic $\Z[\Z^{2n}]$--module
$E_8 \times T^{2n}$,  pass to a large covering
space $p \colon T^{2n} \to T^{2n}$ to obtain a $\delta$--controlled
quadratic $\Z$--module $\wwtilde{E}_8$ over $T^{2n}$ representing an
element of $L_{2n} (T^{2n}; \Z, \epsilon, \delta)$.
It is simple to verify that  $\forget (\wwtilde{E}_8) = p^! \, (E_8 \times T^{2n})$,
where $p^!$ is the $L$--theory transfer.  The transfer invariance
results discussed in \fullref{S-transfer} imply that
$\forget (\wwtilde{E}_8) = E_8 \times T^{2n}$.  Thus,
$\wwtilde{E}_8$ gives a $\delta$--controlled realization of the form
$E_8$ over $T^{2n}$.

\subsection{Controlled surgery obstructions}

\begin{defn}
Let $p\colon X \to B$ be a map to a metric space $B$ and
$\epsilon>0$. A map $f\colon Y \to X$ is an {\em $\epsilon$--homotopy
equivalence} over $B$, if there exist a map $g\colon X \to Y$
and homotopies $H_t$ from $g\circ f$ to $1_Y$ and $K_t$ from
$f\circ g$ to $1_X$  such that $\text{diam}\,(p\circ f
\circ H_t(y))<\epsilon$ for every $y\in Y$, and $\text{diam}
\,(p\circ K_t(x)) <\epsilon$, for every $x\in X$.  This means that
the tracks of $H$ and $K$ are $\epsilon$--small as viewed from $B$.
\end{defn}

Controlled surgery theory addresses the question of the existence and
uniqueness of controlled manifold structures on a space.  Polyhedra
homotopy equivalent to compact topological manifolds satisfy
the Poincar\'e duality isomorphism.  Likewise, there is a notion of
$\epsilon$--Poincar\/e duality satisfied by polyhedra finely
equivalent to a manifold.  Poincar\'e duality can be estimated by the
diameter of cap product with a fundamental class as a chain homotopy
equivalence.

\begin{defn} \label{D-poincare}
Let $p\colon X \to B$ be a map, where $X$ is a polyhedron and
$B$ is a metric space. $X$ is an
{\em $\epsilon$--Poincar\'e complex of formal dimension $n$\/}
over $B$ if there exist a subdivision of $X$ such that simplices have
diameter $\ll \epsilon$ in $B$ and an $n$--cycle $y$ in the
simplicial chains of $X$ so that $\cap y \colon C^{\sharp} (X) \to
C_{n-\sharp} (X)$ is an $\epsilon$--chain homotopy equivalence in the
sense that
$\cap y$ and the chain homotopies have the property that the image of
each generator $\sigma$ only involves generators whose images under
$p$ are within an $\epsilon$--neighborhood of $p(\sigma)$ in $B$.
\end{defn}

To formulate simply-connected controlled surgery problems, the
notion of locally trivial fundamental group from the viewpoint of
the control space is needed.  This can be formalized using the
notion of $UV^1$ maps as follows.

\begin{defn} \label{D-uvone} Given $\delta > 0$,
a map $p\colon X \to B$ is called $\delta$--$UV^1$ if for
any polyhedral pair $(P,Q)$, with $\text{dim}\, (P) \leqslant 2$, and
maps $\alpha_0\colon Q \to X$ and $\beta \colon P \to B$ such
that $p \circ \alpha_0 = \beta |_{Q}$,
\[
\xymatrix{
 Q \ar[r]^{\alpha_0} \ar[d]_i  & X \ar[d]^p \\
P \ar@{.>}[ur]^\alpha \ar[r]_\beta & B
 }
\]
there is a map $\alpha\colon
P \to X$ extending $\alpha_0$ so that $p\circ\alpha$ is
$\delta$--homotopic to $\beta$ over $B$.  The map $p$ is
$UV^1$ if it is $\delta$--$UV^1$, for every $\delta > 0$.
\end{defn}

Let $B$ be a compact metric ENR and $n \geqslant 5$.  Given
$\epsilon > 0$, there is a $\delta > 0$ such that if
$p \colon X \to B$ is a $\delta$--Poincar\'e duality space
over $B$ of formal dimension $n$,
$(f,b) \colon M^n \to X$ is a surgery problem, and
$p$ is $\delta$--$UV^1$. By the Stability \fullref{stability}
there is a well-defined surgery obstruction
$$\sigma_*(f,b) \in \varprojlim_{\epsilon}\varprojlim_{\delta}
L_n(B;\Z,\epsilon,\delta) = H_n (B;\mathbb{L})$$
such that $(f,b)$ is normally cobordant to an $\epsilon$--homotopy
equivalence for any $\epsilon>0$ if and only if $\sigma_*(f,b) = 0$.
See Ranicki and Yamasaki \cite{RY04} for an exposition of controlled
$L$--theory.

The main theorem of Pedersen--Quinn--Ranicki \cite{PQR-zero-one} is the
following controlled surgery exact
sequence (see also Ferry \cite{sF-zero-three} and Ranicki--Yamasaki
\cite{RY04}).

\begin{thm}
Suppose $B$ is a compact metric ENR and $n \geqslant 4$.  There is a
stability threshold $\epsilon_0 > 0$ such that for any
$0 < \epsilon < \epsilon_0$, there is $\delta > 0$ with the property
that if $p \colon N \to B$ is a $\delta$--$UV^1$ map, with $N$ is a
compact $n$--manifold, there is an exact sequence
\[
H_{n+1} (B; \mathbb{L}) \to \mathcal{S}_{\epsilon,\delta} (N,f) \to
[N, \partial N;  G/\mathit{TOP}, \ast] \to H_n (B; \mathbb{L}).
\]
\end{thm}

Here, $\mathcal{S}_{\epsilon,\delta}$ is the controlled structure
set defined as the set of equivalence classes of pairs $(M,g)$,
where $M$ is a topological manifold and $g \colon (M, \partial M)
\to (N, \partial N)$ restricts to a homeomorphism on $\partial N$
and is a $\delta$--homotopy equivalence relative to the boundary.
The pairs $(M_1,g_1)$ and $(M_2, g_2)$ are equivalent if there
is a homeomorphism $h \colon M_1 \to M_2$ such that
$g_1$ and $h \circ g_2$ are $\epsilon$--homotopic rel boundary.
As in classical surgery, the map $H_{n+1}(B; \mathbb{L}) \to
\mathcal{S}_{\epsilon,\delta} (N,f)$ is defined using controlled
Wall realization.

\section{Exotic homology manifolds}\label{S-exotic}

In \cite{BFMW93}, exotic ENR homology manifolds of
dimensions greater than 5 are constructed as limits of sequences
of controlled Poincar\'e complexes $\{X_i, i \geqslant 0\}$.
These complexes are
related by maps $p_i \colon X_{i+1} \to X_i$ such that
$X_{i+1}$ is $\epsilon_{i+1}$--Poincar\'e over $X_i$,
$i \geqslant 0$, and $p_i$ is an
$\epsilon_i$--homotopy equivalence over $X_{i-1}$, $i \geqslant 1$,
where $\sum \epsilon_i < \infty$.  Beginning, say, with a closed
manifold $X_0$, the sequence $\{X_i\}$ is constructed iteratively
using cut-paste constructions on closed manifolds.  The gluing
maps are obtained using the Wall realization of controlled surgery
obstructions, which emerge as a non-trivial local
index in the limiting ENR homology manifold.  As pointed out in
the Introduction, our main goal is to give an explicit construction of
the first controlled stage $X_1$ of this construction using the
quadratic form $E_8$, beginning with the $2n$--dimensional torus
$X_0 = T^{2n}$, $n \geqslant 3$.  The construction of subsequent stages
follows from fairly general arguments presented in \cite{BFMW93}
and leads to an index--9 ENR homology manifold not
homotopy equivalent to any closed topological manifold.
Since an explicit algebraic description of the controlled quadratic
module $\wwtilde{E}_8$ over $T^{2n}$ has already been given in
\fullref{S-E-eight}, we conclude the paper with a review of how
this quadratic module can be used to construct $X_1$.

\begin{figure}[ht!]
\labellist\small
\pinlabel {$C$} [bl] at 90 104
\pinlabel {$D$} at 55 12
\pinlabel {$N$} [tl] at 223 5
\endlabellist
\centerline{\includegraphics[width=1.8in]{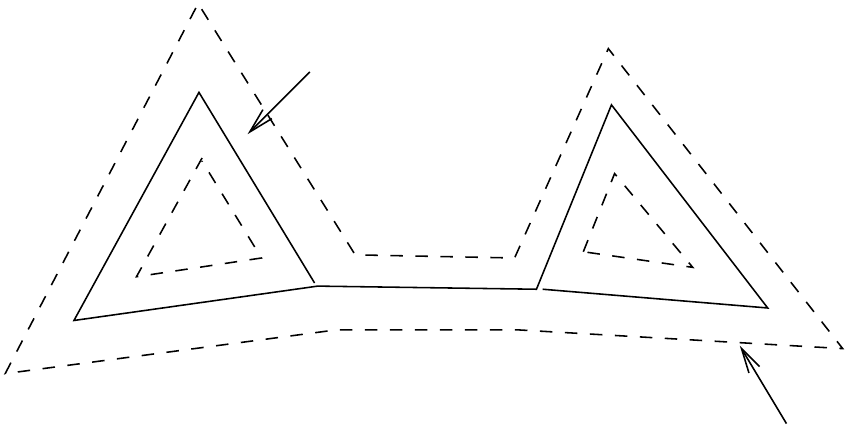}}
\caption{}
\label{F-cut}
\end{figure}

Let $P$ be the 2--skeleton of a fine triangulation of
$T^{2n}$, and $C$ a regular neighborhood of $P$ in
$T^{2n}$.  The closure of the complement of $C$ in $T^{2n}$
will be denoted $D$, and the common boundary
$N = \partial C = \partial D$ (see \fullref{F-cut}).
Given $\delta > 0$, we may assume that the inclusions of
$C, D$ and $N$ into $T^{2n}$ are all $\delta$--$UV^1$
by taking a fine enough triangulation.

Let  $(K, \varphi)$ be a geometric $\Z$--module over
$T^{2n}$ representing the controlled quadratic form $\wwtilde{E}_8$,
where $K \cong \Z^r$ is a free $\Z$--module with basis
$S = \{e_1, \ldots, e_r \}$ and $\varphi \colon S \to T^{2n}$
is a map.  If $Q \subset T^{2n}$ is the dual complex of $P$,
after a small perturbation, we can assume that
$\varphi (S) \cap \left( P \cup Q \right) = \emptyset$.  Composing
this deformation with a retraction $T^{2n} \setminus
\left( P \cup Q \right) \to N$, we can assume that $\varphi$
factors through $N$, that is, the geometric module is actually
realized over $N$.

Using a controlled analogue of the Wall Realization Theorem
\cite[Theorem 5.8]{cW-seventy} applied to the identity map of $N$, realize this
quadratic module over $N \subset T^{2n}$ to obtain
a degree-one normal map $F \colon
(V, N, N') \to (N \times I , N\times \{ 0 \}, N \times
\{ 1 \})$ satisfying:
\begin{itemize}
\item[(a)] $F|_{N} = 1_N$.
\item[(b)] $f = F |_{N'} \colon N' \to N$ is a fine homotopy
equivalence over $T^{2n}$.
\item[(c)] The controlled surgery obstruction of $F$ rel
$\partial$ over $T^{2n}$ is $\wwtilde{E}_8 \!\in\!
H_{2n}(T^{2n};{\mathbb L})$.
\end{itemize}
The map $F$ can be assumed to be $\delta$--$UV^1$
using controlled analogues of $UV^1$ deformation results
of Bestvina and Walsh \cite{kK95}.

\begin{figure}[ht!]
\begin{center}
\labellist\small
\pinlabel {$C$} at 25 25
\pinlabel {$C_f$} at 65 25
\pinlabel {$-V$} at 105 25
\pinlabel {$D$} at 145 25
\endlabellist
\includegraphics[width=1.7in]{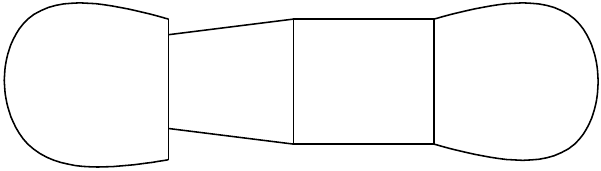}
\end{center}
\caption{The Poincar\'e complex $X_1$}
\label{F-xone}
\end{figure}

Let $C_f$ be the mapping cylinder of $f$.  Form a Poincar\'e
complex $X_1$ by pasting $C_f \cup_{N'} (-V)$ into
$T^{2n}$ along $N$, that is,
\[
X_1 = C \cup_N C_f \cup_{N'} (-V) \cup_N D,
\]
as shown in \fullref{F-xone}.
Our next goal is to define the map $p_1 \colon
X_1 \to X_0 = T^{2n}$.

Let $g \colon N \to N'$ be a controlled homotopy inverse
of $f$.  Composing $f$ and $g$, and using an estimated version of the
Homotopy Extension Theorem (see e.g. \cite{BFMW93})
and the controlled Bestvina--Walsh Theorem, one can
modify $F$ to a $\delta$--$UV^1$ map $G \colon V \to
C_g$, so that $G |_{N'} = 1_{N'}$ and $G |_N = 1_N$.

\begin{figure}[ht!]
\begin{center}
\labellist\small
\pinlabel {$C$} at 25 110
\pinlabel {$C_f$} at 65 110
\pinlabel {$-V$} at 105 110
\pinlabel {$D$} at 145 110
\pinlabel {$\id$} at 35 70
\pinlabel {$\id$} at 75 70
\pinlabel {$G$} at 110 70
\pinlabel {$\id$} at 155 70
\pinlabel {$C$} at 25 25
\pinlabel {$C_f$} at 65 25
\pinlabel {$C_g$} at 105 25
\pinlabel {$D$} at 145 25
\endlabellist
\includegraphics[width=1.7in]{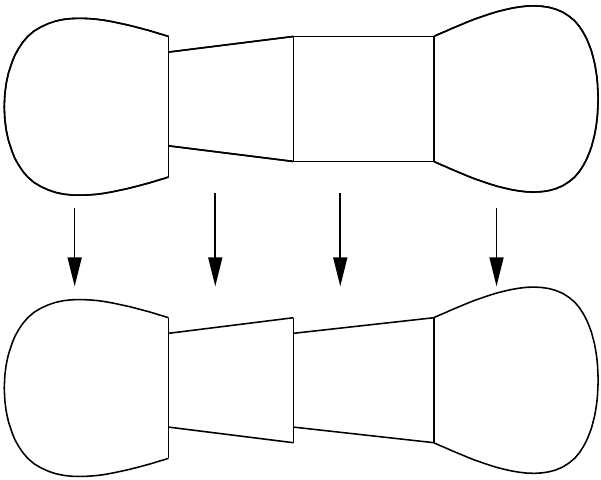}
\end{center}
\caption{The map $p^\ast_1 \colon X_1 \to X_1'$}
\label{F-map}
\end{figure}

Let $X_1' = C \cup_N C_f \cup_{N'} C_g \cup_N D$ and
$p_1^\ast \colon X_1 \to X_1'$ be as indicated in
\fullref{F-map}.
Crushing $C_f \cup_{N'} C_g$ to
$N = \partial C$, we obtain the desired map
$p_1 \colon X_1 \to T^{2n} = C \cup_N D$.

\begin{figure}[t]
\labellist\small
\pinlabel {$C$} at 25 135
\pinlabel {$V$} at 60 135
\pinlabel {$N'\!{\times}I$} at 91 135
\pinlabel {$-V$} at 120 135
\pinlabel {$D$} at 160 135
\pinlabel {$\id$} at 30 85
\pinlabel {$F$} at 65 85
\pinlabel {$\pi$} at 95 85
\pinlabel {$\id$} at 135 85
\pinlabel {$\id$} at 165 85
\pinlabel {$C$} at 25 30
\pinlabel {$N\!{\times}I$} at 61.5 30
\pinlabel {$C_f$} at 90 30
\pinlabel {$-V$} at 120 30
\pinlabel {$D$} at 160 30
\endlabellist
\centerline{\includegraphics[width=2in]{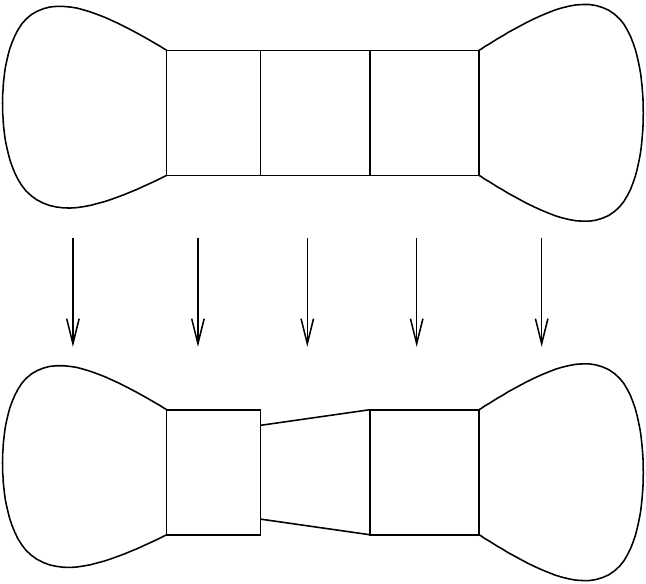}}
\caption{The map $\phi \colon M \to X_1$}
\label{F-index}
\end{figure}

To conclude, as in \cite{BFMW96}, we argue  that
$X_1$ is not homotopy equivalent to any closed
topological manifold.  To see this, consider the closed
manifold
\[
M = C \cup_N V \cup_{N'} N' \times I \cup_{N'}
(-V) \cup_N D
\]
and the degree-one normal map $\phi \colon M \to X_1$
depicted in \fullref{F-index}, where $\pi \colon N' \times I
\to C_f$ is induced by $f \colon N' \to N$.
The controlled
surgery obstruction of $\phi$ over $T^{2n}$ is the generator
\begin{multline*}
\sigma_*(\phi) = E_8 \times T^{2n}=(0,\dots,0,1)
\in L_0(\Z)= L_{2n}(\Z[\Z^{2n}])^\mathit{INV}\\
\subset H_{2n}\left(T^{2n};\LL\right)=L_{2n}(\Z[\Z^{2n}])=\sum\limits^{2n}_{r=0}
\di{2n \choose r}L_{2n-r}(\Z)
\end{multline*}
of the subgroup of the transfer invariant elements (\fullref{invariance}).
Let $\LL\langle 1 \rangle$ be the 1--connective cover of $\LL$,
the simply-connected surgery spectrum with 0th space (homotopy equivalent to)
$G/\mathit{TOP}$. Now
$$L_{2n}(\Z[\Z^{2n}]) = H_{2n}(T^{2n};\LL)
= H_{2n}(T^{2n};\LL\langle 1 \rangle ) \oplus L_0(\Z)$$\nopagebreak
with
$$H_{2n}(T^{2n};\LL\langle 1 \rangle) = [T^{2n},G/\mathit{TOP}] = 
\sum\limits^{2n}_{r=1}\di{2n \choose r}L_{2n-r}(\Z) \subset
L_{2n}(\Z[\Z^{2n}])$$ the subgroup of the surgery obstructions of
normal maps $M_1 \to T^{2n}$. The surgery obstruction of any
normal map $\phi_1\co M_1 \to X_1$ is of the type
$$\sigma_*(\phi_1) = (\tau,1) \neq 0 \in L_{2n}(\Z[\Z^{2n}]) = 
[T^{2n},G/\mathit{TOP}] \oplus L_0(\Z)$$
for some $\tau \in [T^{2n},G/\mathit{TOP}]$, since the variation of
normal invariant only changes the component of the surgery
obstruction in $[T^{2n},G/\mathit{TOP}] \subset  L_{2n}(\Z[\Z^{2n}])$.
Thus, $X_1$ is not homotopy equivalent to any topological
manifold. In the terminology of Ranicki \cite[Chapter 17]{ranicki3} the
total surgery obstruction $s(X_1)\in {\mathscr S}_{2n}(X_1)$ has
image
$$(p_1)_*s(X_1)  = 1 \in {\mathscr S}_{2n}(T^{2n}) = L_0(\Z).$$
The Bryant--Ferry--Mio--Weinberger procedure for constructing an ENR
homology manifold starting with $p_1 \colon X_1 \to T^{2n}$ leads
to a homology manifold homotopy equivalent to $X_1$.  Thus, from
the quadratic form $E_8$, we obtained a compact  index--9 ENR
homology $2n$--manifold $\mathfrak{X}_8$ which is not homotopy
equivalent to any closed topological manifold.

\subsection*{Acknowledgment}
This research was partially supported by NSF grant DMS-0071693.

\bibliographystyle{gtart}
\bibliography{link}

\begin{thebibliography}{}
\providecommand\bibmarginpar{\leavevmode\marginpar}
\def\urlstyle#1{{\tt #1}}

\bibitem{BFMW93}
\textbf{J Bryant}, \textbf{S Ferry}, \textbf{W Mio}, \textbf{S Weinberger},
  \href{http://dx.doi.org/10.1090/S0273-0979-1993-00381-1} {\emph{Topology of
  homology manifolds}}, Bull. Amer. Math. Soc. $($N.S.$)$ 28 (1993) 324--328
  \xox{MR}{1183997}

\bibitem{BFMW96}
\textbf{J Bryant}, \textbf{S Ferry}, \textbf{W Mio}, \textbf{S Weinberger},
  \href{http://dx.doi.org/10.2307/2118532} {\emph{Topology of homology
  manifolds}}, Ann. of Math. $(2)$ 143 (1996) 435--467 \xox{MR}{1394965}

\bibitem{jC78}
\textbf{J\,W Cannon}, \href{http://dx.doi.org/10.1090/S0002-9904-1978-14527-3}
  {\emph{The recognition problem: what is a topological manifold?}}, Bull.
  Amer. Math. Soc. 84 (1978) 832--866 \xox{MR}{0494113}

\bibitem{jC79}
\textbf{J\,W Cannon}, \href{http://dx.doi.org/10.2307/1971245} {\emph{Shrinking
  cell-like decompositions of manifolds. {C}odimension three}}, Ann. of Math.
  $(2)$ 110 (1979) 83--112 \xox{MR}{541330}

\bibitem{clauwens2}
\textbf{F Clauwens}, \emph{The $K$--theory of almost symmetric forms}, from:
  ``Topological structures, II (Proc. Sympos. Topology and Geom., Amsterdam,
  1978), Part 1'', Math. Centre Tracts 115, Math. Centrum, Amsterdam (1979)
  41--49 \xox{MR}{565824}

\bibitem{clauwens3}
\textbf{F Clauwens}, \emph{Product formulae for surgery obstructions}, from:
  ``Algebraic topology, Aarhus 1978'', Lecture Notes in Math. 763, Springer,
  Berlin (1979)  198--211 \xox{MR}{561223}

\bibitem{clauwens-one}
\textbf{F Clauwens}, \emph{Surgery on products I, {II}}, Nederl. Akad.
  Wetensch. Indag. Math. 41 (1979) 121--132, 133--144 \xox{MR}{535561}

\bibitem{rD86}
\textbf{R\,J Daverman}, \emph{Decompositions of manifolds}, Pure and Applied
  Mathematics 124, Academic Press, Orlando, FL (1986) \xox{MR}{872468}

\bibitem{rE80}
\textbf{R\,D Edwards}, \emph{The topology of manifolds and cell-like maps},
  from: ``Proceedings of the International Congress of Mathematicians
  (Helsinki, 1978)'', Acad. Sci. Fennica, Helsinki (1980)  111--127
  \xox{MR}{562601}

\bibitem{sF-zero-three}
\textbf{S\,C Ferry}, \emph{Epsilon-Delta surgery over $\mathbb{Z}$}, preprint,
  Rutgers University (2003)

\bibitem{G}
\textbf{M Gromov}, \emph{Positive curvature, macroscopic dimension, spectral
  gaps and higher signatures}, from: ``Functional analysis on the eve of the
  21st century, Vol.\ II (New Brunswick, NJ, 1993)'', Progr. Math. 132,
  Birkh\"auser, Boston (1996)  1--213 \xox{MR}{1389019}

\bibitem{kK95}
\textbf{K Kawamura}, \href{http://dx.doi.org/10.1016/0166-8641(94)00038-5}
  {\emph{An inverse system approach to {M}enger manifolds}}, Topology Appl. 61
  (1995) 281--292 \xox{MR}{1317082}

\bibitem{MR}
\textbf{R\,J Milgram}, \textbf{A\,A Ranicki}, \emph{The $L$--theory of
  {L}aurent extensions and genus 0 function fields}, J. Reine Angew. Math. 406
  (1990) 121--166 \xox{MR}{1048238}

\bibitem{wM00}
\textbf{W Mio}, \emph{Homology manifolds}, from: ``Surveys on surgery theory,
  Vol. 1'', Ann. of Math. Stud. 145, Princeton Univ. Press, Princeton, NJ
  (2000)  323--343 \xox{MR}{1747540}

\bibitem{PQR-zero-one}
\textbf{E\,K Pedersen}, \textbf{F Quinn}, \textbf{A Ranicki}, \emph{Controlled
  surgery with trivial local fundamental groups}, from: ``High-dimensional
  manifold topology'', World Sci. Publishing, River Edge, NJ (2003)  421--426
  \xox{MR}{2048731}

\bibitem{PY}
\textbf{E Pedersen}, \textbf{M Yamasaki},
  \href{http://dx.doi.org/10.2140/gtm.2006.9.69} {\emph{Stability in controlled
  $L$--theory}}, from: ``Exotic homology manifolds (Oberwolfach 2003)'', Geom.
  Topol. Monogr. 9 (2006)  69--88 \xox{arXiv}{math.GT/0402218}

\bibitem{fQ83}
\textbf{F Quinn}, \href{http://dx.doi.org/10.1007/BF01389323}
  {\emph{Resolutions of homology manifolds, and the topological
  characterization of manifolds}}, Invent. Math. 72 (1983) 267--284
  \xox{MR}{700771}

\bibitem{fQ87}
\textbf{F Quinn}, \href{http://dx.doi.org/10.1307/mmj/1029003559} {\emph{An
  obstruction to the resolution of homology manifolds}}, Michigan Math. J. 34
  (1987) 285--291 \xox{MR}{894878}

\bibitem{ranicki-one}
\textbf{A Ranicki}, \href{http://dx.doi.org/10.1112/plms/s3-40.1.87} {\emph{The
  algebraic theory of surgery I: {F}oundations}}, Proc. London Math. Soc. $(3)$
  40 (1980) 87--192 \xox{MR}{560997}

\bibitem{ranicki3}
\textbf{A\,A Ranicki}, \emph{Algebraic $L$--theory and topological manifolds},
  Cambridge Tracts in Mathematics 102, Cambridge University Press, Cambridge
  (1992) \xox{MR}{1211640}

\bibitem{ranicki-two}
\textbf{A Ranicki}, \emph{Lower $K$-- and $L$--theory}, London Mathematical
  Society Lecture Note Series 178, Cambridge University Press, Cambridge (1992)
  \xox{MR}{1208729}

\bibitem{ranicki-four}
\textbf{A Ranicki}, \emph{High-dimensional knot theory}, Springer Monographs in
  Mathematics, Springer, New York (1998) \xox{MR}{1713074}

\bibitem{RY04}
\textbf{A Ranicki}, \textbf{M Yamasaki},
  \href{http://dx.doi.org/10.2140/gtm.2006.9.107} {\emph{Controlled
  $L$--theory}}, from: ``Exotic homology manifolds (Oberwolfach 2003)'', Geom.
  Topol. Monogr. 9 (2006)  107--156 \xox{arXiv}{math.GT/0402217}

\bibitem{Sh}
\textbf{J\,L Shaneson}, \href{http://dx.doi.org/10.2307/1970726} {\emph{Wall's
  surgery obstruction groups for $G\times Z$}}, Ann. of Math. $(2)$ 90 (1969)
  296--334 \xox{MR}{0246310}

\bibitem{cW-seventy}
\textbf{C\,T\,C Wall}, \emph{Surgery on compact manifolds}, Mathematical
  Surveys and Monographs 69, American Mathematical Society, Providence, RI
  (1999) \xox{MR}{1687388}

\bibitem{sW95}
\textbf{S Weinberger}, \emph{Nonlocally linear manifolds and orbifolds}, from:
  ``Proceedings of the International Congress of Mathematicians, Vol.\ 1, 2
  (Z\"urich, 1994)'', Birkh\"auser, Basel (1995)  637--647 \xox{MR}{1403964}

\bibitem{mY87}
\textbf{M Yamasaki}, \href{http://dx.doi.org/10.1007/BF01391832}
  {\emph{$L$--groups of crystallographic groups}}, Invent. Math. 88 (1987)
  571--602 \xox{MR}{884801}

\end{thebibliography}

\end{document}